# On the Iwasawa Invariants of Elliptic Curves

Ralph Greenberg
University of Washington
and
Vinayak Vatsal
University of British Columbia

## 1 Introduction

Suppose that $E$ is an elliptic curve defined over $\mathbb{Q}$, and that $p$ is a prime where $E$ has good ordinary reduction. Under the assumption that $E$ is modular, one can define nonnegative integers $\lambda_E^{\text{alg}}, \mu_E^{\text{alg}}, \lambda_E^{\text{anal}}$, and $\mu_E^{\text{anal}}$. The "algebraic" Iwasawa invariants $\lambda_E^{\text{alg}}$ and $\mu_E^{\text{alg}}$ are defined in terms of the structure of the $p$-primary subgroup $\text{Sel}_E(\mathbb{Q}_\infty)_p$ of the Selmer group for $E$ over the cyclotomic $\mathbb{Z}_p$-extension $\mathbb{Q}_\infty$ of $\mathbb{Q}$. The definition of the "analytic" invariants $\mu_E^{\text{anal}}$ and $\lambda_E^{\text{anal}}$ is in terms of the $p$-adic L-function for $E$ constructed by Mazur and Swinnerton-Dyer [MSD74]. We will recall these definitions below. Our purpose in this article is to prove in certain cases that $\mu_E^{\text{alg}} = \mu_E^{\text{anal}} = 0$ and that $\lambda_E^{\text{alg}} = \lambda_E^{\text{anal}}$. These equalities, together with a deep theorem of Kato, imply the main conjecture for $E$ over $\mathbb{Q}_\infty$. In this introduction, we will discuss the nature of our results and give an outline of the proofs for the case of modular elliptic curves. We want to point out, however, that the theorems proven in the text apply to modular forms, as well as to elliptic curves with multiplicative reduction at $p$. We will not attempt to state the most general versions here.

If $K$ is any algebraic extension of $\mathbb{Q}$, then the Selmer group for $E$ over $K$ is a certain subgroup of $H^1(G_K, E(\overline{\mathbb{Q}})_{\text{tors}})$, where $G_K = \text{Gal}(\overline{\mathbb{Q}}/K)$. The Selmer group fits into an exact sequence

$$0 \to E(K) \otimes \mathbb{Q}/\mathbb{Z} \to \text{Sel}_E(K) \to \Sha_E(K)$$

where $\Sha_E(K)$ denotes the Shafarevich-Tate group for $E$ over $K$. Let $K = \mathbb{Q}_\infty$. Then $\Gamma = \text{Gal}(\mathbb{Q}_\infty/\mathbb{Q})$ acts on $\text{Sel}_E(\mathbb{Q}_\infty)$. Its $p$-primary subgroup $\text{Sel}_E(\mathbb{Q}_\infty)_p$ (to which we give the discrete topology) can be regarded as a $\Lambda$-module, where $\Lambda = \mathbb{Z}_p[[\Gamma]]$ is the completed group algebra for $\Gamma$ over $\mathbb{Z}_p$. Kato has proven that $\text{Sel}_E(\mathbb{Q}_\infty)_p$ is $\Lambda$-cotorsion, as Mazur conjectured in [Maz72]. That is, the Pontrjagin dual $X_E(\mathbb{Q}_\infty)$ of $\mathrm{Sel}_E(\mathbb{Q}_\infty)_p$ is a torsion $\Lambda$-module. It is an easy result that $X_E(\mathbb{Q}_\infty)$ is finitely generated as a $\Lambda$-module, and so the classification of finitely generated $\Lambda$-modules asserts that one has a pseudo-isomorphism

$$X_E(\mathbb{Q}_\infty) \sim \left(\bigoplus_{i=1}^{n} \Lambda/(f_i(T)^{a_i})\right) \bigoplus \left(\bigoplus_{j=1}^{m} \Lambda/(p^{\mu_j})\right),$$

where one identifies $\Lambda$ with the formal power series ring $\mathbb{Z}_p[[T]]$ in the usual way: $T = \gamma - 1$, where $\gamma$ is a fixed topological generator of $\Gamma$. The $f_i(T)$'s are irreducible distinguished polynomials in $\Lambda$. The $a_i$'s and the $\mu_j$'s are positive integers. One can then define the algebraic Iwasawa invariants by

$$\lambda_E^{\mathrm{alg}} = \sum_{i=1}^{n} a_i \deg(f_i(T)), \text{ and } \mu_E^{\mathrm{alg}} = \sum_{j=0}^{m} \mu_j. \qquad (1)$$

For the formulation of the Main Conjecture, it is also necessary to define the "characteristic polynomial" for the $\Lambda$-module $X_E(\mathbb{Q}_\infty)$; it is given by

$$f_E^{\mathrm{alg}}(T) = p^{\mu_E^{\mathrm{alg}}} \cdot \prod_{i=1}^{n} f_i(T)^{a_i}. \qquad (2)$$

Thus $p^{\mu_E^{\mathrm{alg}}}$ is the exact power of $p$ dividing $f_E^{\mathrm{alg}}(T)$ in $\Lambda$, and $\lambda_E^{\mathrm{alg}} = \deg(f_E^{\mathrm{alg}}(T))$.

One can also describe $\lambda_E^{\mathrm{alg}}$ in the following group-theoretic way:

$$\left(\mathrm{Sel}_E(\mathbb{Q}_\infty)_p\right)_{\mathrm{div}} \cong (\mathbb{Q}_p/\mathbb{Z}_p)^{\lambda_E^{\mathrm{alg}}}.$$

To deduce this fact from the above definitions, one uses the fact that $\Lambda/(f(T))$ is a free $\mathbb{Z}_p$-module of rank equal to $\deg(f(T))$, for any distinguished polynomial $f(T)$. The roots of $f_E^{\mathrm{alg}}(T)$ (counting multiplicities) are the eigenvalues of $\gamma - 1 = T$ acting on the vector space $X_E(\mathbb{Q}_\infty) \otimes_{\mathbb{Z}_p} \overline{\mathbb{Q}}_p$, which has dimension $\lambda_E^{\mathrm{alg}}$. The invariant $\lambda_E^{\mathrm{alg}}$ can be quite large, as some of our later examples and remarks will show. On the other hand, it is expected that $\mu_E^{\mathrm{alg}}$ is usually zero. But this is *not* always the case. Indeed, Mazur showed in [Maz72] that $\mu_E^{\mathrm{alg}}$ can be positive for certain $E$ and $p$. Here are some of the known results concerning $\mu_E^{\mathrm{alg}}$:

**A.** Suppose that $E_1$ and $E_2$ are elliptic curves defined over $\mathbb{Q}$. Let $p$ be an odd prime where $E_1$ and $E_2$ have good ordinary reduction. Assume that $E_1[p] \cong E_2[p]$ as Galois modules. Then $\mathrm{Sel}_{E_1}(\mathbb{Q}_\infty)[p]$ is finite if and only if $\mathrm{Sel}_{E_2}(\mathbb{Q}_\infty)[p]$ is finite. Consequently, if $\mathrm{Sel}_{E_1}(\mathbb{Q}_\infty)_p$ is $\Lambda$-cotorsion and $\mu_{E_1}^{\mathrm{alg}} = 0$, then $\mathrm{Sel}_{E_2}(\mathbb{Q}_\infty)_p$ is $\Lambda$-cotorsion and $\mu_{E_2}^{\mathrm{alg}} = 0$.

**B.** Suppose that $E$ is an elliptic curve over $\mathbb{Q}$, and that $p$ is an odd prime of good ordinary reduction. Assume that $E$ admits a cyclic $\mathbb{Q}$-isogeny of degree $p^t$ with kernel $\Phi$. Assume that the action of $G_{\mathbb{Q}}$ on $\Phi$ is ramified at $p$ and odd. (That is, the action of an inertia group $I_p$ is nontrivial, and the action of complex conjugation is by $-1$.) Then $\mu_E^{\text{alg}} \geq t$.

**C.** Suppose that $E$ is an elliptic curve defined over $\mathbb{Q}$ with good ordinary reduction at an odd prime $p$. Assume that $E$ admits a $\mathbb{Q}$-isogeny of degree $p$ with kernel $\Phi$, and that the action of $G_{\mathbb{Q}}$ on $\Phi$ is either ramified at $p$ and even, or unramified at $p$ and odd. Then $\mu_E^{\text{alg}} = 0$.

The first result is rather easy to prove, and can be generalized to elliptic curves over any number field $K$. Let $K_\infty = K\mathbb{Q}_\infty$, the cyclotomic $\mathbb{Z}_p$-extension of $K$. If $E$ is an elliptic curve over $K$, then Mazur conjectured that $\text{Sel}_E(K_\infty)$ should be $\Lambda$-cotorsion. The only really general result in this direction is due to Kato, who showed that $\text{Sel}_E(K_\infty)_p$ is $\Lambda$-cotorsion if $E$ is defined over $\mathbb{Q}$ and is modular and $K/\mathbb{Q}$ is an abelian extension. The second result can be proven by using results of Poitou and Tate on the local and global Euler characteristics for the Galois cohomology of $\Phi$. It can also be generalized (with a rather more complicated statement) to any number field $K$. The proof of the third result is based on the well-known theorem of Ferrero and Washington on the vanishing of the classical $\mu$-invariant for abelian extensions of $\mathbb{Q}$. In this case one can prove that the $\mu$-invariant of $\text{Sel}_E(K_\infty)$ vanishes if $E$ satisfies the hypotheses of **C** and if $K$ is a totally real abelian number field. One can find a rather thorough discussion of these results in [Gre99], where the case $p = 2$ is also treated. Also, the above results are valid if $E$ has multiplicative reduction at $p$.

We now define the analytic invariants. Suppose that $E$ is a modular elliptic curve over $\mathbb{Q}$, and that $p$ is a prime of good ordinary reduction. For any Dirichlet character $\rho$, we let $L(E/\mathbb{Q}, \rho, s)$ denote the the Hasse-Weil L-series for $E$, twisted by the character $\rho$. Let $\Omega_E$ denote the real Néron period for $E$. If $\rho$ is even, then it is known that $L(E/\mathbb{Q}, \rho, 1)/\Omega_E \in \overline{\mathbb{Q}}$, where we have fixed an embedding of $\overline{\mathbb{Q}}$ into $\mathbb{C}$. We also fix an embedding of $\overline{\mathbb{Q}}$ into $\overline{\mathbb{Q}}_p$. Mazur and Swinnerton-Dyer have constructed an element $\mathscr{L}(E/\mathbb{Q}, T) \in \Lambda \otimes \mathbb{Q}_p$ satisfying a certain interpolation property which we now describe. Suppose that $\rho \in \text{Hom}(\Gamma, \mu_{p^\infty})$ is a character of finite order. Since $\gamma$ is a topological generator of $\Gamma$, $\rho$ is determined by $\rho(\gamma) = \zeta \in \mu_{p^\infty}$. One can view $\rho$ as a Dirichlet character of $p$-power order and conductor. Assuming that $\rho$ is nontrivial, its conductor is of the form $p^m$, and, assuming that $p$ is odd, $\zeta$ has order $p^{m-1}$. Then $\mathscr{L}(E/\mathbb{Q}, T)$ is characterized by

$$\mathscr{L}(E/\mathbb{Q}, \zeta - 1) = \tau(\rho^{-1}) \cdot \alpha_p^{-m} \cdot \frac{L(E/\mathbb{Q}, \rho, 1)}{\Omega_E}. \tag{3}$$

where $\rho$ runs over the set of all nontrivial characters of $\Gamma$. Here $\tau(\rho^{-1})$ denotes the usual Gauss sum, and $\alpha_p$ denotes the eigenvalue for Frobenius acting on the maximal (1-dimensional) unramified quotient of the $p$-adic Tate module of $E$. Alternatively, $\alpha_p$ is the unit solution of the equations $\alpha_p \beta_p = p$ and $\alpha_p + \beta_p = 1 + p - \#\widetilde{E}(\mathbb{F}_p)$, where $\mathbb{F}_p$ is the finite field with $p$ elements and $\widetilde{E}$ is the reduction of $E$ at $p$. The power series $\mathscr{L}(E/\mathbb{Q}, T)$ is determined by this property because a nonzero element of $\Lambda \otimes \mathbb{Q}_p$ has only finitely many zeroes. Now a theorem of Rohrlich [Roh89] states that $L(E/\mathbb{Q}, \rho, 1) \neq 0$ for all but finitely many $\rho$ of $p$-power conductor. Thus the element $\mathscr{L}(E/\mathbb{Q}, T)$ is nonzero. Using the Weierstrass preparation theorem, we define the invariants $\mu_E^{\text{anal}}$ and $\lambda_E^{\text{anal}}$ by writing

$$\mathscr{L}(E/\mathbb{Q}, T) = p^{\mu_E^{\text{anal}}} \cdot u(T) \cdot f(T)$$

where $f(T)$ is a distinguished polynomial of degree $\lambda_E^{\text{anal}}$, and $u(T)$ is an invertible power series. We define $f_E^{\text{anal}}(T) = p^{\mu_E^{\text{anal}}} \cdot f(T)$, where $f(T)$ is as above. One should have $\mu_p^{\text{anal}} \geq 0$. That is, $f_E^{\text{anal}}(T)$ should be in $\mathbb{Z}_p[T]$. This is known when $E[p]$ is irreducible as a Galois module and also when $E[p]$ is reducible under certain additional hypotheses. The integrality of $\mathscr{L}(E/\mathbb{Q}, T)$ is discussed in Theorem 3.7.

If $p$ is odd, one can identify $\Gamma$ with $\text{Gal}\,(\mathbb{Q}(\mu_{p^\infty})/\mathbb{Q}(\mu_p))$. Let $\chi$ denote the cyclotomic character giving the action of $\text{Gal}\,(\mathbb{Q}(\mu_{p^\infty})/\mathbb{Q})$ on $\mu_{p^\infty}$. We let $\kappa = \chi|_\Gamma$, which induces an isomorphism $\Gamma \cong 1 + p\mathbb{Z}_p$. Then the $p$-adic L-function $L_p(E/\mathbb{Q}, s)$ is defined by

$$L_p(E/\mathbb{Q}, s) = \mathscr{L}(E/\mathbb{Q}, \kappa(\gamma)^{1-s} - 1).$$

Even though the power series $\mathscr{L}(E/\mathbb{Q}, T)$ depends on the choice of $\gamma$, the function $L_p(E/\mathbb{Q}, s)$ is independent of this choice. Also, $L_p(E/\mathbb{Q}, 1) = \mathscr{L}(E/\mathbb{Q}, 0)$, which was not specified above. If $E$ has good reduction at $p$, it turns out that one has

$$L_p(E/\mathbb{Q}, 1) = \mathscr{L}(E/\mathbb{Q}, 0) = (1 - \alpha_p^{-1})^2 \frac{L(E/\mathbb{Q}, 1)}{\Omega_E},$$

where $\alpha_p$ is as above. Note that $(1 - \alpha_p^{-1}) = (1 - \beta_p p^{-1})$ is part of the Euler factor for $p$ in $L(E/\mathbb{Q}, s)$, evaluated at $s = 1$. Now the Main Conjecture can be stated as follows:

**Conjecture (1.1) (Mazur)** *We have $f_E^{alg}(T) = f_E^{anal}(T)$*

Obviously, this would imply that $\mu_E^{\text{alg}} = \mu_E^{\text{anal}}$, and that $\lambda_E^{\text{alg}} = \lambda_E^{\text{anal}}$. Kato has proven the following weaker statement:

**Theorem (1.2) (Kato)** *The polynomial $f_E^{alg}(T)$ divides $f_E^{anal}(T)$ in $\mathbb{Q}_p[T]$.*

As a consequence, it is clear that the equality $\lambda_E^{\text{alg}} = \lambda_E^{\text{anal}}$ implies that $f_E^{\text{alg}}(T)$ and $f_E^{\text{anal}}(T)$ differ by multiplication by a power of $p$. The further equality $\mu_E^{\text{alg}} = \mu_E^{\text{anal}}$ then implies the Main Conjecture.

A proof of Kato's theorem in the form just stated (for modular elliptic curves with good, ordinary reduction at a prime $p > 2$) has been presented by Scholl and Rubin. (See [Sch98] and [Rub98].) A more general version is contained in [Kat] which applies also to cusp forms where $p$ divides the level and of arbitrary weight. Under certain hypotheses, Kato even proves the divisibility in $\mathbb{Z}_p[T]$.

Our results may be stated as follows.

**Theorem (1.3)** *Assume that $E$ is a modular elliptic curve over $\mathbb{Q}$, and that $p$ is an odd prime where $E$ has good ordinary reduction. Assume also that $E$ admits a $\mathbb{Q}$-isogeny of degree $p$ with kernel $\Phi$, and that the action of $G_\mathbb{Q}$ on $\Phi$ is either ramified at $p$ and even, or unramified at $p$ and odd. Then $\lambda_E^{alg} = \lambda_E^{anal}$ and $\mu_E^{alg} = \mu_E^{anal} = 0$.*

**Theorem (1.4)** *Assume that $E_1$ and $E_2$ are modular elliptic curves over $\mathbb{Q}$. Let $p$ be an odd prime where $E_1$ and $E_2$ have good ordinary reduction. Assume also that $E_1[p] \cong E_2[p]$ as Galois modules, and that these are irreducible. If the equalities $\mu_{E_1}^{alg} = \mu_{E_1}^{anal} = 0$ and $\lambda_{E_1}^{alg} = \lambda_{E_1}^{anal}$ hold, then so do the equalities $\mu_{E_2}^{alg} = \mu_{E_2}^{anal} = 0$ and $\lambda_{E_2}^{alg} = \lambda_{E_2}^{anal}$.*

The set of primes where these theorems apply is rather limited. The hypotheses in theorem (1.3) can hold only when $p = 3, 5, 7, 13$, or $37$. However, for the first four of these primes, there are infinitely many distinct $j$-invariants $j_E$ which can occur. The hypotheses are also preserved by even quadratic twists of conductor prime to $p$. If one is willing to consider more general modular forms, then one can prove a similar theorem for any $p$ dividing certain Bernoulli numbers.

As an illustrative example of theorem (1.3), take $p = 5$, and let $J$ denote any of the three nonisomorphic elliptic curves of conductor 11. It is known that the $G_\mathbb{Q}$-module $J[5]$ is reducible and has composition factors isomorphic to $\mu_5$ and $\mathbb{Z}/5\mathbb{Z}$ as $G_\mathbb{Q}$-modules. Let $\psi$ be an odd quadratic character such that $\psi(5) \neq 0$, corresponding to the imaginary quadratic field $\mathbb{Q}(\sqrt{-c})$, and let $E = J_{-c}$ denote the associated quadratic twist. The curve $E$ has a cyclic $\mathbb{Q}$-isogeny of degree 5 with kernel isomorphic to either $\mu_5 \otimes \psi$ or $(\mathbb{Z}/5\mathbb{Z}) \otimes \psi$, so that theorem (1.3) clearly applies. Thus $\mu_E^{alg} = \mu_E^{anal} = 0$, and, as we will later prove, $\lambda_E^{alg} = \lambda_E^{anal} = 2\lambda_\psi + \epsilon_\psi$, where $\lambda_\psi$ denotes the classical $\lambda$-invariant for the imaginary quadratic field $\mathbb{Q}(\sqrt{-c})$ and for the prime $p = 5$, and $\epsilon_\psi = 0$ if 11 is inert or ramified in this field, $\epsilon_\psi = 1$ if 11 splits. There are examples (found by T. Fukuda) where $\lambda_\psi$ is quite large. The record so far is $\lambda_\psi = 10$ for $c = 3,624,233$. In this case $\lambda_{anal} = \lambda_{alg} = 21$. As we have remarked above, the Main Conjecture is valid for $E$. Unfortunately, if $E$ is a quadratic twist of $J$ by an even character, we can prove very little. If $J = X_0(11)$, then both $\mu_E^{alg}$ and $\mu_E^{anal}$ are positive. If $J$ is the curve 11C in Cremona [Cre92], then we would expect $\mu_E^{alg}$ and $\mu_E^{anal}$ to be equal to zero. We can prove neither equality. Even if this were true, we don't know how to prove the equality $\lambda_E^{alg} = \lambda_E^{anal}$.

As an example for theorem (1.4), consider the two elliptic curves

$$E_1 : y^2 = x^3 + x - 10, \qquad E_2 : y^2 = x^3 - 584x + 5444.$$

Then $E_1$ and $E_2$ are modular, and have conductors 52 and 364 respectively. We will take $p = 5$. We find that $E_1[p] \cong E_2[p]$ by comparing the $q$-expansions of the associated modular forms, and by observing that neither curve admits a rational 5-isogeny. Both curves have good, ordinary reduction at 5. Now one has $L(E_1/\mathbb{Q}, 1)/\Omega_{E_1} = 1/2 \in \mathbb{Z}_5^\times$. Also, $\#\widetilde{E}_1(\mathbb{F}_5) = 4$, and $\alpha_5 \equiv 2 \pmod{5\mathbb{Z}_5}$. Thus $1 - \beta_5 5^{-1} = 1 - \alpha_5^{-1} \not\equiv 0 \pmod{5\mathbb{Z}_5}$, so that $\mathscr{L}(E_1/\mathbb{Q}, T)$ must be in $\Lambda^\times$. This means that $f_{E_1}^{\text{anal}}(T) = 1$ and $\mu_{E_1}^{\text{anal}} = \lambda_{E_1}^{\text{anal}} = 0$. On the other hand, a well-known result of Kolyvagin shows that $\text{Sel}_{E_1}(\mathbb{Q})_5 = 0$. One can deduce from this fact (together with the facts that $\#E_1(\mathbb{F}_5) = 4$, and that the Tamagawa factors $c_2 = 1, c_3 = 2$, are all 5-adic units) that $\text{Sel}_{E_1}(\mathbb{Q}_\infty)^\Gamma = 0$. This implies that $\text{Sel}_{E_1}(\mathbb{Q}_\infty) = 0$, and so $f_{E_1}^{\text{alg}}(T) = 1$, $\lambda_{E_1}^{\text{alg}} = \mu_{E_1}^{\text{alg}} = 0$. Thus the hypotheses of theorem (1.4) are satisfied, and we can conclude that $\mu_{E_2}^{\text{alg}} = \mu_{E_2}^{\text{anal}} = 0$, and that $\lambda_{E_2}^{\text{alg}} = \lambda_{E_2}^{\text{anal}}$. The proof of Theorem (1.4) will show that $\lambda_{E_2}^{\text{alg}}$ and $\lambda_{E_2}^{\text{anal}}$ are both equal to 5. More precisely, it turns out that we have

$$0 \to E_2(\mathbb{Q}_\infty) \otimes \mathbb{Q}_5/\mathbb{Z}_5 \to \text{Sel}_{E_2}(\mathbb{Q}_\infty)_5 \to \text{III}_{E_2}(\mathbb{Q}_\infty)_5 \to 0$$

with $E_2(\mathbb{Q}_\infty) \otimes \mathbb{Q}_5/\mathbb{Z}_5 = \mathbb{Q}_5/\mathbb{Z}_5$ and $\text{III}_{E_2}(\mathbb{Q}_\infty)_5 = (\mathbb{Q}_5/\mathbb{Z}_5)^4$. The Mordell-Weil group $E_2(\mathbb{Q})$ has rank one, and the first group above is just the image of $E_2(\mathbb{Q}) \otimes \mathbb{Q}_5/\mathbb{Z}_5$ under the restriction map $\text{Sel}_{E_2}(\mathbb{Q})_5 \to \text{Sel}_{E_2}(\mathbb{Q}_\infty)_5$.

In [RS93], Rubin and Silverberg have shown that there are infinitely many elliptic curves $E/\mathbb{Q}$ such that $E[3]$ or $E[5]$ has a given structure as $G_\mathbb{Q}$-modules (assuming that at least one such curve exists). They describe the family of such curves in terms of a rational parametrization. In particular, the curves $E$ such that $E[5] \cong E_1[5]$, where $E_1$ is as above, are given by

$$E : y^2 = x^3 + a(t)x + b(t),$$

where $a(t)$ and $b(t)$ are explicitly specified polynomials over $\mathbb{Q}$ (of degree 20 and 30 respectively) and $t \in \mathbb{Q}$ satisfies $\Delta(t) = 4a(t)^3 + 27b(t)^2 \neq 0$. One can then verify that each of these curves is modular and has good ordinary reduction at $p = 5$. Theorem (1.4) then implies that $\mu_E^{\text{anal}} = \mu_E^{\text{alg}} = 0$, and that $\lambda_E^{\text{anal}} = \lambda_E^{\text{alg}}$. The Main Conjecture also holds for each such $E$ and $p = 5$. Interestingly, it can be shown that the $\lambda$-invariant is unbounded in the family.

We will now outline the structure of the proofs. Let $E$ be an elliptic curve over $\mathbb{Q}$ and let $p$ be a prime of good ordinary reduction for $E$. The Selmer group $\text{Sel}_E(\mathbb{Q}_\infty)_p$ is defined as the kernel of the following "global-to-local" map

$$H^1(\mathbb{Q}_\Sigma/\mathbb{Q}_\infty, E[p^\infty]) \to \prod_{\ell \in \Sigma} \mathcal{H}_\ell(\mathbb{Q}_\infty), \tag{4}$$

where $\Sigma$ denotes any finite set of primes containing $p, \infty$, and the primes of bad reduction for $E$. For each finite prime $\ell \in \Sigma$, the group $\mathcal{H}_\ell(\mathbb{Q}_\infty)$ is defined by

$$\mathcal{H}_\ell(\mathbb{Q}_\infty) = \prod_{\eta | \ell} H^1((\mathbb{Q}_\infty)_\eta, E[p^\infty])/\mathrm{im}\,(\kappa_\eta).$$

Here $\eta$ runs over the finite set of places of $\mathbb{Q}_\infty$ over $\ell$, and $(\mathbb{Q}_\infty)_\eta$ denotes the union of the completions of the finite layers $\mathbb{Q}_n$ at $\eta$. Also, $\kappa_\eta$ denotes the local Kummer map

$$\kappa_\eta : E((\mathbb{Q}_\infty)_\eta) \otimes \mathbb{Q}_p/\mathbb{Z}_p \to H^1((\mathbb{Q}_\infty)_\eta, E[p^\infty]).$$

If $\ell = \infty$ and $p$ is odd, then we simply take $\mathcal{H}_\ell(\mathbb{Q}_\infty) = 0$. (A more careful definition would be necessary if $p = 2$, especially because $\infty$ splits completely in $\mathbb{Q}_\infty/\mathbb{Q}$. We will not go into this here.) For each $\ell \in \Sigma$, the group $\mathcal{H}_\ell(\mathbb{Q}_\infty)$ is a cofinitely generated $\Lambda$-module. The group $\mathcal{H}_p(\mathbb{Q}_\infty)$ has $\Lambda$-corank 1. For $\ell \neq p$, it can be shown that $\mathcal{H}_\ell(\mathbb{Q}_\infty)$ is $\Lambda$-cotorsion, and has $\mu$-invariant 0. Its structure is not hard to study. A helpful and easily verified fact is that $E((\mathbb{Q}_\infty)_\eta) \otimes \mathbb{Q}_p/\mathbb{Z}_p = 0$ when $\eta$ is prime to $p$, so that the image of $\kappa_\eta$ is zero. Hence $\mathcal{H}_\ell(\mathbb{Q}_\infty)$ is simply the product $\prod_{\eta | \ell} H^1((\mathbb{Q}_\infty)_\eta, E[p^\infty])$. It turns out in fact that $\mathcal{H}_\ell(\mathbb{Q}_\infty) \cong (\mathbb{Q}_p/\mathbb{Z}_p)^{\sigma_E^{(\ell)}}$, where $\sigma_E^{(\ell)}$ is a non-negative integer which is easily determined from the Euler factor for $\ell$ in $L(E/\mathbb{Q}, s)$. Let $h_E^{(\ell)}(T) \in \Lambda$ denote the characteristic polynomial for the $\Lambda$-module $\mathcal{H}_\ell(\mathbb{Q}_\infty)\widehat{\ }$, where the hat indicates the Pontrjagin dual. This polynomial has degree $\sigma_E^{(\ell)}$.

Let $\Sigma_0$ be any subset of $\Sigma$ which does not contain $p$. We define the corresponding "non-primitive" Selmer group by

$$\mathrm{Sel}_E^{\Sigma_0}(\mathbb{Q}_\infty)_p = \ker\left(H^1(\mathbb{Q}_\Sigma/\mathbb{Q}_\infty, E[p^\infty]) \to \prod_{\ell \in \Sigma - \Sigma_0} \mathcal{H}_\ell(\mathbb{Q}_\infty)\right).$$

This group will play a crucial role in our arguments. Its main advantage is that if $\Sigma_0$ is chosen to contain all primes of bad reduction, and if $p$ is odd, then $\mathrm{Sel}_E^{\Sigma_0}(\mathbb{Q}_\infty)_p[p]$ is determined completely by the $G_\mathbb{Q}$-module $E[p]$. We will be more precise about this later. Assuming that $\mathrm{Sel}_E(\mathbb{Q}_\infty)_p$ is $\Lambda$-cotorsion (which has been proven by Kato when $E$ is modular and has good or multiplicative reduction at $p$), one can show that the map (4) defining the Selmer group is surjective. It follows that

$$\mathrm{Sel}_E^{\Sigma_0}(\mathbb{Q}_\infty)_p / \mathrm{Sel}_E(\mathbb{Q}_\infty)_p \cong \prod_{\ell \in \Sigma_0} \mathcal{H}_\ell(\mathbb{Q}_\infty). \tag{5}$$

Consequently, if we denote the characteristic polynomial of the $\Lambda$-module $\mathrm{Sel}_E^{\Sigma_0}(\mathbb{Q}_\infty)_p\widehat{\ }$ by $f_{E,\Sigma_0}^{\mathrm{alg}}(T)$, then we have

$$f_{E,\Sigma_0}^{\mathrm{alg}}(T) = f_E^{\mathrm{alg}}(T) \prod_{\ell \in \Sigma_0} h_E^{(\ell)}(T), \tag{6}$$

and the $\lambda$-invariant of $\text{Sel}_E^{\Sigma_0}(\mathbb{Q}_\infty)_p$, which we denote by $\lambda_{E,\Sigma_0}^{\text{alg}}$, is given by

$$\lambda_{E,\Sigma_0}^{\text{alg}} = \lambda_E^{\text{alg}} + \sum_{\ell \in \Sigma_0} \sigma_E^{(\ell)}. \tag{7}$$

As for the $\mu$-invariant, it is obvious that $\mu_E^{\text{alg}} = \mu_{E,\Sigma_0}^{\text{alg}}$. If $\mu_E^{\text{alg}} = 0$, then, as we will explain later, it turns out that $\text{Sel}_E^{\Sigma_0}(\mathbb{Q}_\infty)_p$ is a divisible group if $\Sigma_0$ is chosen as above. In particular,

$$\text{Sel}_E^{\Sigma_0}(\mathbb{Q}_\infty)_p \cong (\mathbb{Q}_p/\mathbb{Z}_p)^{\lambda_{E,\Sigma_0}^{\text{alg}}},$$

and therefore we find that $\lambda_{E,\Sigma_0}^{\text{alg}} = \dim_{\mathbb{F}_p}(\text{Sel}_E^{\Sigma_0}(\mathbb{Q}_\infty)_p[p])$. As we mentioned above, $\text{Sel}_E^{\Sigma_0}(\mathbb{Q}_\infty)[p]$ is determined by the Galois module $E[p]$ if $\Sigma_0$ is chosen suitably. Thus $E[p]$ determines $\lambda_{E,\Sigma_0}^{\text{alg}}$ and we can then recover $\lambda_E^{\text{alg}}$ from (7). Note that $E[p]$ does not determine $\lambda_E^{\text{alg}}$. This is clear from the example of $E_1$ and $E_2$ given previously.

There is a similar story for the $p$-adic L-functions $L_p(E/\mathbb{Q}, s)$. We assume that $E$ is modular, and that it has good ordinary reduction at $p$. Let $\Sigma_0$ again be any finite set of primes not containing $p$. For any Dirichlet character $\rho$, we denote by $L^{\Sigma_0}(E/\mathbb{Q}, \rho, s)$ the nonprimitive complex L-function formed from $L(E/\mathbb{Q}, \rho, s)$ by simply omitting the Euler factors for the primes in $\Sigma_0$. If $\ell$ is any prime, we denote the corresponding Euler factor by $P_\ell(E/\mathbb{Q}, \rho, \ell^{-s})$, where

$$P_\ell(E/\mathbb{Q}, \rho, X) = (1 - \rho(\ell)\alpha_\ell X)(1 - \rho(\ell)\beta_\ell X),$$

for the usual quantities $\alpha_\ell$ and $\beta_\ell$. (Possibly one or both of $\alpha_l$, $\beta_l$ are zero.) Then

$$L^{\Sigma_0}(E/\mathbb{Q}, \rho, s) = L(E/\mathbb{Q}, \rho, s) \prod_{\ell \in \Sigma_0} P_\ell(E/\mathbb{Q}, \rho, \ell^{-s}).$$

It is easy to modify $\mathscr{L}(E/\mathbb{Q}, T)$ to construct a nonprimitive $p$-adic L-function. We want to define an element $\mathscr{L}^{\Sigma_0}(E/\mathbb{Q}, T) \in \Lambda \otimes \mathbb{Q}_p$ by requiring that

$$\mathscr{L}^{\Sigma_0}(E/\mathbb{Q}, \zeta - 1) = \tau(\rho^{-1}) \cdot \alpha_p^{-m} \cdot \frac{L^{\Sigma_0}(E/\mathbb{Q}, \rho, 1)}{\Omega_E}$$

for each nontrivial character $\rho \in \text{Hom}(\Gamma, \mu_{p^\infty})$, with the notation as in (3). For $\ell \in \Sigma_0$, let $\gamma_\ell$ denote the Frobenius automorphism for $\ell$ in $\Gamma = \text{Gal}(\mathbb{Q}_\infty/\mathbb{Q})$. (Note that $\ell$ is unramified in $\mathbb{Q}_\infty/\mathbb{Q}$, since $\ell \neq p$.) Consider the element

$$\mathscr{P}_\ell = (1 - \alpha_\ell \ell^{-1} \gamma_\ell)(1 - \beta_\ell \ell^{-1} \gamma_\ell) \in \mathbb{Z}_p[[\Gamma]].$$

We will write the element $\mathscr{P}_\ell$ of $\mathbb{Z}_p[[\Gamma]]$ as a power series $\mathscr{P}_\ell(T)$. To do this we replace $\gamma_\ell$ by $(1+T)^{f_\ell}$, where $f_\ell \in \mathbb{Z}_p$ is determined by $\gamma^{f_\ell} = \gamma_\ell$. Since $\rho(\gamma_\ell) = \rho(\ell)$ when $\rho$ is viewed as a Dirichlet character, it follows that

$$\mathscr{P}_\ell(\zeta - 1) = P_\ell(E/\mathbb{Q}, \rho, \ell^{-1}).$$

Thus we simply define

$$\mathscr{L}^{\Sigma_0}(E/\mathbb{Q}, T) = \mathscr{L}(E/\mathbb{Q}, T) \prod_{\ell \in \Sigma_0} \mathscr{P}_\ell(T).$$

Note that $\mathscr{P}_\ell(T)$ is a nonzero element of $\Lambda$ which is not divisible by $p$. We will show later that $\mathscr{P}_\ell(T)$ generates the characteristic ideal of the module $\mathcal{H}_\ell(\mathbb{Q}_\infty)\widehat{\phantom{x}}$. If we define the polynomial $f^{\mathrm{anal}}_{E,\Sigma_0}(T)$ in the same way as in the primitive case, it follows that

$$f^{\mathrm{anal}}_{E,\Sigma_0}(T) = f^{\mathrm{anal}}_E(T) \prod_{\ell \in \Sigma_0} h^{(\ell)}_E(T). \tag{8}$$

The degree $\lambda^{\mathrm{anal}}_{E,\Sigma_0}$ of $f^{\mathrm{anal}}_{E,\Sigma_0}(T)$ is given by

$$\lambda^{\mathrm{anal}}_{E,\Sigma_0} = \lambda^{\mathrm{anal}}_E + \sum_{\ell \in \Sigma_0} \sigma^{(\ell)}_E. \tag{9}$$

One can also define the $\mu$-invariant $\mu^{\mathrm{anal}}_{E,\Sigma_0}$ in the obvious way. One clearly has $\mu^{\mathrm{anal}}_{E,\Sigma_0} = \mu^{\mathrm{anal}}_E$.

Combining the above observations with the previous considerations yields the following result.

**Theorem (1.5)** *Assume that $E$ is a modular elliptic curve with good ordinary reduction at $p$. Let $\Sigma_0$ be any finite set of nonarchimedean primes not containing $p$. Then the following equivalences hold:*

1. $\mu^{alg}_E = \mu^{anal}_E \iff \mu^{alg}_{E,\Sigma_0} = \mu^{anal}_{E,\Sigma_0}$

2. $\lambda^{alg}_E = \lambda^{anal}_E \iff \lambda^{alg}_{E,\Sigma_0} = \lambda^{anal}_{E,\Sigma_0}$.

3. $f^{alg}_E(T) = f^{anal}_E(T) \iff f^{alg}_{E,\Sigma_0}(T) = f^{anal}_{E,\Sigma_0}(T)$.

This theorem is crucial to our arguments because it is only the nonprimitive $p$-adic $L$-functions and the nonprimitive Selmer groups corresponding to a suitable choice of the set $\Sigma_0$ which behave well under congruences. The observation that a "main conjecture" should be equivalent to a nonprimitive analogue is an old one. In the context of classical Iwasawa theory (involving Kubota-Leopoldt $p$-adic $L$-functions) such an equivalence was proved in [Gre77]

Now we can describe the proof of Theorem (1.4). Suppose that $E_1$ and $E_2$ satisfy the stated conditions. By theorem (1.5), we have $\lambda^{\text{alg}}_{E_1,\Sigma_0} = \lambda^{\text{anal}}_{E_1,\Sigma_0}$. Choose $\Sigma$ to be a finite set of primes containing $p, \infty$, and all primes where either $E_1$ or $E_2$ has bad reduction. Let $\Sigma_0 = \Sigma - \{p, \infty\}$. We are assuming that $p$ is odd and that $E_1[p] \cong E_2[p]$. Since $\mu^{\text{alg}}_{E_1} = 0$, it follows that $\mu^{\text{alg}}_{E_1,\Sigma_0} = 0$, and that $\text{Sel}^{\Sigma_0}_{E_1}(\mathbb{Q}_\infty)_p[p]$ is finite. Thus $\text{Sel}^{\Sigma_0}_{E_2}(\mathbb{Q}_\infty)[p]$ is finite too, which implies that $\mu^{\text{alg}}_{E_2,\Sigma_0} = \mu^{\text{alg}}_{E_2} = 0$. We obtain the equality of $\lambda^{\text{alg}}_{E_1,\Sigma_0}$ and $\lambda^{\text{alg}}_{E_2,\Sigma_0}$, since, as remarked earlier, these now depend only on $E_1[p] \cong E_2[p]$. The analogous results for the analytic invariants are proved by a theory of congruences for L-functions. More precisely, we show that

$$\mathscr{L}^{\Sigma_0}(E_1/\mathbb{Q}, T) \equiv u \cdot \mathscr{L}^{\Sigma_0}(E_2/\mathbb{Q}, T) \pmod{p\Lambda}, \tag{10}$$

where $u \in \mathbb{Z}_p^\times$. Since the vanishing of $\mu^{\text{anal}}_{E_1}$ implies the vanishing of $\mu^{\text{anal}}_{E_1,\Sigma_0}$, we see that $\mu^{\text{anal}}_{E_2,\Sigma_0}$ and $\mu^{\text{anal}}_{E_2}$ vanish as well. It is then clear that the degrees of $f^{\text{anal}}_{E_1,\Sigma_0}(T)$ and $f^{\text{anal}}_{E_2,\Sigma_0}(T)$ coincide. Here we use the fact that if $\mathscr{L}(T) \in \Lambda$ is any power series such that $p \nmid \mathscr{L}(T)$, then $\mathscr{L}(T) = u(T) \cdot f(T)$ where $f(T)$ is a distinguished polynomial, $u(T)$ is in $\Lambda^\times$, and the degree of $f(T)$ is determined by the image of $\mathscr{L}(T)$ in $\Lambda/p\Lambda$. Thus we obtain the following equalities:

$$\lambda^{\text{alg}}_{E_2,\Sigma_0} = \lambda^{\text{alg}}_{E_1,\Sigma_0} = \lambda^{\text{anal}}_{E_1,\Sigma_0} = \lambda^{\text{anal}}_{E_2,\Sigma_0}.$$

Then by theorem (1.5) we obtain $\lambda^{\text{alg}}_{E_2} = \lambda^{\text{anal}}_{E_2}$ as stated in theorem (1.4).

Now assume that $E$ satisfies the hypotheses of Theorem (1.3). Then there is an exact sequence

$$0 \to \Phi \to E[p] \to \Psi \to 0 \tag{11}$$

of $G_\mathbb{Q}$-modules, where $\Phi$ and $\Psi$ are cyclic of order $p$. Let $\varphi$ and $\psi$ denote the corresponding $\mathbb{F}_p^\times$-valued characters of $G_\mathbb{Q}$. We may view these as taking values in $\mathbb{Z}_p^\times$, and then we have $\varphi\psi = \omega$, where $\omega$ is the usual Teichmüller character. Replacing $E$ by an isogenous elliptic curve if necessary (which turns out not to affect any of the Iwasawa invariants), we may assume that $\varphi$ is ramified at $p$ and even, so that $\psi$ is unramified at $p$ and odd. Consider the $G_\mathbb{Q}$-modules $C = \mu_{p^\infty} \otimes \psi^{-1}$ and $D = (\mathbb{Q}_p/\mathbb{Z}_p) \otimes \psi$. Each is isomorphic to $\mathbb{Q}_p/\mathbb{Z}_p$ as a group, and we have $\Phi = C[p]$ and $\Psi = D[p]$. Classical Iwasawa theory, which involves the study of certain natural Galois groups regarded as $\Lambda$-modules, can be reformulated in terms of certain subgroups $S_C(\mathbb{Q}_\infty)$ and $S_D(\mathbb{Q}_\infty)$ (which we also refer to as "Selmer groups") of the Galois cohomology groups $H^1(\mathbb{Q}_\infty, C)$ and $H^1(\mathbb{Q}_\infty, D)$ respectively. The classical Main Conjecture (which was proven by Mazur and Wiles) gives a precise connection between the structure of these Selmer groups and the $p$-adic L-functions attached to $C$ and $D$. These $p$-adic $L$-functions are essentially the Kubota-Leopoldt $p$-adic L-functions, and are defined by elements $\mathscr{L}(C,T)$ and $\mathscr{L}(D,T)$ of $\Lambda$ satisfying a certain interpolation property. These power series are closely related: we have $\mathscr{L}(C,T) = \mathscr{L}(D,T^\iota)$, where $T^\iota = (1+T)^{-1} - 1$. (The ring $\Lambda$ has a natural involution $\iota$ induced by $\gamma \to \gamma^{-1}$.) The Ferrero-Washington theorem asserts that the $\mu$-invariants of these power series as well as those of the $\Lambda$-modules $S_C(\mathbb{Q}_\infty)\hat{\ }$ and $S_D(\mathbb{Q}_\infty)\hat{\ }$ are all equal to zero. The corresponding analytic $\lambda$-invariants are equal: $\lambda_C^{\text{anal}} = \lambda_D^{\text{anal}}$. Furthermore, one has $S_C(\mathbb{Q}_\infty) \cong S_D(\mathbb{Q}_\infty)^\iota$ as $\Lambda$-modules, where the $\iota$ indicates changing of the $\Lambda$-module structure by the involution $\iota$. Thus we get $\lambda_C^{\text{alg}} = \lambda_D^{\text{alg}}$. The main conjecture for $C$ and $D$ shows that $\lambda_C^{\text{anal}} = \lambda_C^{\text{alg}} = \lambda_D^{\text{alg}} = \lambda_D^{\text{anal}}$. The analogue of theorem (1.5) is also valid. Let $\Sigma_0$ be any set of primes not containing $p$. One can define non-primitive Selmer groups as before, as well as non-primitive $p$-adic L-functions. The $\mu$-invariants will still be zero, and the equality $\lambda_C^{\text{anal}} = \lambda_C^{\text{alg}}$ implies $\lambda_{C,\Sigma_0}^{\text{anal}} = \lambda_{C,\Sigma_0}^{\text{alg}}$. The corresponding statement also holds for $D$. However, one does not usually have $\lambda_{C,\Sigma_0} = \lambda_{D,\Sigma_0}$.

Suppose now that $\Sigma$ is a finite set of primes containing $p, \infty$, and the primes of bad reduction for E. Let $\Sigma_0 = \Sigma - \{p, \infty\}$. The exact sequence relating $E[p]$ to $\Phi$ and $\Psi$ will then allow us to prove that $\mu_E^{\text{alg}} = 0$ (as a consequence of the Ferrero-Washington theorem) and that $\lambda_{E,\Sigma_0}^{\text{alg}} = \lambda_{C,\Sigma_0}^{\text{alg}} + \lambda_{D,\Sigma_0}^{\text{alg}}$. As for the analytic side, we will use (11) to produce the congruence

$$\mathscr{L}^{\Sigma_0}(E/\mathbb{Q},T) \equiv u\mathscr{L}^{\Sigma_0}(C,T) \cdot \mathscr{L}^{\Sigma_0}(D,T) \pmod{p\Lambda}, \tag{12}$$

where $u \in \mathbb{Z}_p^\times$. The Ferrero-Washington theorem again shows that $\mu_E^{\text{anal}} = \mu_{E,\Sigma_0}^{\text{anal}} = 0$ and the above congruence then implies that $\lambda_{E,\Sigma_0}^{\text{anal}} = \lambda_{C,\Sigma_0}^{\text{anal}} + \lambda_{D,\Sigma_0}^{\text{anal}}$. We conclude from the Main Conjecture for $C$ and $D$ that $\lambda_{E,\Sigma_0}^{\text{alg}} = \lambda_{E,\Sigma_0}^{\text{anal}}$. The equality $\lambda_E^{\text{alg}} = \lambda_E^{\text{anal}}$ is then a consequence of theorem (1.5).

We want to briefly describe the technique used to produce congruences between $p$-adic L-functions as in (10) and (12). The basic idea goes back to fundamental work of Mazur (see [Maz77] and [Maz79]), where it was made clear that congruences between analytic L-values could be studied using the Hecke-module structure of the cohomology of modular curves. It turns out that what one needs is a certain "multiplicity-one" result for this cohomology, and this was subsequently provided by work of Mazur, Ribet, Wiles, and others. The result (10) was essentially proven in [Vat97], where multiplicity-one was used to define canonical periods associated to modular forms. The result (12) has its origins in the paper [Maz79], where the case $J_0(N)$ was treated, when $N$ is a prime. In particular, Mazur considered the curve $E = X_0(11)$. His key idea was that the exact sequence (11) may be viewed as giving a congruence between the weight 2 modular form associated to $E$ and a certain Eisenstein series $G(z)$ whose associated Dirichlet series $L(G, s)$ is given by $L(G, s) = L(\psi, s) \cdot L(\varphi \omega^{-1}, s - 1)$. The congruence (12) is the translation of this fact into the setting of L-functions: the L-function of $E$ is congruent to the L-function of $G(z)$. This idea of Mazur was subsequently generalized by Stevens to arbitrary level in [Ste82], but Stevens' result was subject to certain unverified hypotheses (see [Ste82], page 109, and Sec. 4.2.) These restrictions were replaced by similar (but easily verified) conditions in [Vat97]. For the purposes of this paper, we have to go still further and analyze the "canonical periods" associated to modular forms in [Vat97]. When the modular form $f$ corresponds to an elliptic curve $E$, we show that the canonical periods of $f$ are essentially the Néron periods of a certain curve (the *optimal* curve) in the isogeny class of $E$.

Fundamental to this approach is the fact that the congruences (10) and (12) involve the *nonprimitive* L-functions $L^{\Sigma_0}(\cdot, T)$. This is the analytic counterpart of the fact that the Galois module $E[p]$ determines only the nonprimitive Selmer group. Consider the first congruence (10). Thus fix two elliptic curves $E_i/\mathbb{Q}, i = 1, 2$, together with an isomorphism $E_1[p] \cong E_2[p]$ of (irreducible) Galois modules. Let $f_i(z)$ denote the newform of level $N_i$ corresponding to $E_i$. We assume that $(N_i, p) = 1$, and that the $E_i$ have ordinary reduction at $p$. Let $\Sigma_0$ denote the set of primes dividing $N_1 N_2$. Our hypothesis $E_1[p] \cong E_2[p]$ implies that if $f_i(z) = \sum a_i(n)q^n$, then we have the congruence $a_1(n) \equiv a_2(n) \pmod{p}$ for all $n$ with $(n, N_1 N_2) = 1$. Let

$$g_i(z) = \sum_{(n, N_1 N_2) = 1} a_i(n) q^n = \sum_{n \geq 1} b_i(n) q^n.$$

We then have $L(g_i, s) = \Sigma b_i(n) n^{-s} = L^{\Sigma_0}(E_i/\mathbb{Q}, s)$. Furthermore, there exists an integer $M$, divisible only by the primes in $\Sigma_0$, such that each $g_i$ is a modular form of level $M$ and weight 2, an eigenform for all the Hecke operators $T_l, U_q$. Hence there are ring homomorphisms $\epsilon_i : \mathbf{T} \to \mathbb{Z}$ such that $\epsilon_i(T_l) = b_i(l)$ for $l \notin \Sigma_0$ and $\epsilon_i(U_q) = 0$ for all $q \in \Sigma_0$, where $\mathbf{T} = \mathbb{Z}[T_l, U_q]$ is the Hecke algebra of level $M$, weight 2. We have $\epsilon_1(T) \equiv \epsilon_2(T) \pmod{p}$ for all $T \in \mathbf{T}$ and therefore there is a unique maximal ideal $\mathfrak{m}$ of $\mathbf{T}$ containing $\ker(\epsilon_1)$, $\ker(\epsilon_2)$ and $p$. The residue field $\mathbf{T}/\mathfrak{m}$ is just $\mathbb{F}_p$. Thus $T_l \equiv b_1(l) \equiv b_2(l) \pmod{\mathfrak{m}}$ for all $l \notin \Sigma_0$ and $U_q \equiv 0 \pmod{\mathfrak{m}}$ for $q \in \Sigma_0$.

Let $\mathbf{T}_\mathfrak{m}$ denote the completion of $\mathbf{T}$ at $\mathfrak{m}$, which is a direct factor in the semilocal $\mathbb{Z}_p$-algebra $\mathbf{T} \otimes_\mathbb{Z} \mathbb{Z}_p$. Thus $\mathbf{T}_\mathfrak{m}$ is a free $\mathbb{Z}_p$-module of finite rank and the hypothesis that $E_i[p]$ is irreducible implies that $\mathbf{T}_\mathfrak{m}$ is Gorenstein. (This is corollary 2 of theorem 2.1 in [Wil95]. Earlier versions of the Gorenstein property for Hecke rings were proved by Mazur, Ribet, Edixhoven and others. Wiles' version is the most general. Equivalently, $H^1(X_1(M), \mathbb{Z})_\mathfrak{m}^\pm$ is a free $\mathbf{T}_\mathfrak{m}$-module of rank 1. Here $H^1(X_1(M), \mathbb{Z})$ is the usual Betti cohomology of the modular curve $X_1(M)$ with coefficients in $\mathbb{Z}$. The superscript $\pm$ indicates the eigenspaces for the action of complex conjugation, and the subscript $\mathfrak{m}$ denotes the $\mathfrak{m}$-adic completion with respect to the natural structure as a $\mathbf{T}$-module.

Let $S_2(\Gamma_1(M), \mathbb{Z})$ denote the cusp forms of weight 2 for $\Gamma_1(M)$ with Fourier coefficients in $\mathbb{Z}$. It is known that $S_2(\Gamma_1(M), \mathbb{Z})$ can be identified with $\mathrm{Hom}(\mathbf{T}, \mathbb{Z})$ as $\mathbf{T}$-modules (DI95, proposition 12.4.1). Thus, $S_2(\Gamma_1(M), \mathbb{Z})_\mathfrak{m} \cong \mathrm{Hom}(\mathbf{T}_\mathfrak{m}, \mathbb{Z}_p)$ and the Gorenstein property of $\mathbf{T}_\mathfrak{m}$ states that $\mathrm{Hom}(\mathbf{T}_m, \mathbb{Z}_p)$ is isomorphic to $\mathbf{T}_m$. Hence $S_2(\Gamma_1(M), \mathbb{Z})_m$ is also a free $\mathbf{T}_m$-module of rank 1 and so there are isomorphisms

$$\theta^\pm : S_2(\Gamma_1(M), \mathbb{Z})_m \xrightarrow{\sim} H^1(X_1(M), \mathbb{Z})_m^\pm,$$

which are equivariant for the action of $\mathbf{T}_m$.

To prove the congruence (10), it suffices to show that

$$\mathscr{L}^{\Sigma_0}(E_1/\mathbb{Q}, \zeta - 1) \equiv u\mathscr{L}^{\Sigma_0}(E_2/\mathbb{Q}, \zeta - 1) \pmod{p\mathbb{Z}_p[\zeta]}$$

for all $p$-power roots of unity $\zeta$ and for some $u \in \mathbb{Z}_\rho^\times$ independent of $\zeta$. This is clear because if $\mathscr{G}(T) = \Sigma c_j T^j \in \Lambda$ is not divisble by $p$ and if $c_\lambda$ is the first unit coefficient, then $\mathscr{G}(\zeta - 1)$ and $(\zeta - 1)^\lambda$ have the same valuation when $\zeta$ has sufficiently large order. Applying this to $\mathscr{G}(T) = \mathscr{L}^{\Sigma_0}(E_1/\mathbb{Q}, T) - u\mathscr{L}^{\Sigma_0}(E_2/\mathbb{Q}, T)$, it follows that $\mathscr{G}(T) \equiv 0 \pmod{p\Lambda}$. Now, as is well-known, one can obtain the values $\tau(\rho^{-1})L(E_i/\mathbb{Q}, \rho, 1)$ from the differential forms $\omega_i(z)^\pm = (g_i(z)dz)^\pm \in H^1(X_1(M), \mathbb{C})^\pm$ (by integration along paths joining the cusps of $X_1(M)(\mathbb{C})$ and forming certain linear combinations). Since $E_1[p] \cong E_2[p]$ implies that $\alpha_p(E_1) \equiv \alpha_p(E_2) \pmod{p\mathbb{Z}_p}$, it is enough to prove a congruence for values of integrals involving $\omega_i(z)^\pm/\Omega_{E_i}^\pm \in H^1(X_1(M), \mathbb{Q})^\pm$, where $\Omega_{E_i}^\pm$ denotes the real or imaginary Néron period for $E_i$. On the other hand, we can define cocycles $\delta_i^\pm = \theta^\pm(g_i) \in H^1(X_1(M), \mathbb{Z})_\mathfrak{m}^\pm$ and we clearly have $\delta_1^\pm \equiv \delta_2^\pm \pmod{p}$. (This means that we have a congruence modulo $pH^1(X_1(M), \mathbb{Z})_\mathfrak{m}^\pm$. We use similar notation elsewhere.) Now $\delta_i^\pm$ must be a rational multiple of $\omega_i(z)^\pm/\Omega_{E_i}^\pm$. To establish the congruences involving $\mathscr{L}(E_i/\mathbb{Q}, \zeta - 1)$, the key result that we need is that there exist $p$-adic units $u_i^\pm$ such that

$$\delta_i^\pm = \frac{\omega_i(z)^\pm}{u_i^\pm \Omega_{E_i}^\pm}.$$

This is the main result of section 3.

There is an analogous explanation for the congruence (12). But the argument is slightly more delicate in this case, as the Gorenstein property is unknown in general. To fix notation, let $E/\mathbb{Q}$ denote an elliptic curve, with corresponding newform $f = \sum a_n q^n$ of level $N$. Assume that there exists an exact sequence of $\mathrm{Gal}(\overline{\mathbb{Q}}/\mathbb{Q})$-modules

$$0 \to \Phi \to E[p] \to \Psi \to 0$$

where $\text{Gal}(\overline{\mathbb{Q}}/\mathbb{Q})$ acts on $\Phi$ and $\Psi$ by abelian characters $\varphi$ and $\psi$ respectively. We assume that the character $\psi$ is odd and unramified. (The case where $\psi$ is even and ramified is similar). Let $\Sigma_0$ denote the set of primes dividing $N$, and let $g = \sum_{(n,N)=1} a_n q^n = \sum b_n q^n$. Then $g$ is an eigenform for all the Hecke operators of a suitable level $M$, and we have $L(g,s) = L_{\Sigma_0}(f,s)$. We let $\mathfrak{m}$ denote the maximal ideal of characteristic $p$ determined by $g$ in the Hecke ring $\mathbf{T}$ of level $M$. Under these hypotheses, it is known that there exists an isomorphism
$$\theta : S_2(\Gamma_1(M), \mathbb{Z}_p)_{\mathfrak{m}} \cong H^1(X_1(M), \mathbb{Z}_p)_{\mathfrak{m}}^+.$$
We define a canonical cocycle $\delta_g = \theta(g) \in H^1(X_1(M), \mathbb{Z}_p)^+$. As before, it can be shown that we have
$$\delta_g = \frac{(g(z)dz)^+}{u \cdot \Omega_E^+}$$
where $\Omega_E^+$ denotes the real Néron period for $E$, and $u$ is a $p$-adic unit. We remark here that although $E$ belongs to a nontrivial class of $p$-isogenous curves, the cocycle $\delta_g$ depends only on $g$, and is independent of the choice of $E$ in the isogeny class. The unit $u$ may depend on the choice of $E$. Implicit in the equality displayed above is the statement that the period $\Omega_E$ is an invariant of the isogeny class, up to $p$-adic unit.

Now let $G(z) = \sum c_n q^n$ denote the Eisenstein series with L-function given by $\sum c_n n^{-s} = L_{\Sigma_0}(\psi, s) \cdot L_{\Sigma_0}(\psi^{-1}, s-1)$. Then we have the congruence $c_n \equiv b_n \pmod{p}$, for $n \geq 1$. Mazur and Stevens have shown how to construct an Eisenstein cocycle $\overline{\delta}_G \in H^1(X_1(M), \mathbb{F}_p)^+[\mathfrak{m}]$, associated to $G$. Here $\mathbb{F}_p$ denotes the finite field with $p$ elements and $H^1(X_1(M), \mathbb{F}_p)^+[\mathfrak{m}]$ denotes the kernel of $\mathfrak{m}$ acting on the Betti cohomology of $X_1(M)$, with coefficients in $\mathbb{F}_p$ (see [Maz79] and [Ste82]). The proof that the Eisenstein cocycle $\overline{\delta}_G$ is nonzero is based on a deep theorem of Washington [Was78].

On the other hand, if we write $\overline{\delta}_g$ for the reduction of the cocycle $\delta_g$, then clearly we have $\overline{\delta}_g \in H^1(X_1(M), \mathbb{F}_p)^+[\mathfrak{m}]$. It turns out that $\overline{\delta}_g$ is nonzero. Now, the $q$-expansion principle implies that $H^1(X_1(M), \mathbb{F}_p)^+[\mathfrak{m}]$ is a 1-dimensional vector space over $\mathbb{F}_p$. This is the crucial part of Mazur's original argument, for it implies that the cocycles $\overline{\delta}_g$ and $\overline{\delta}_G$ coincide up to a factor in $\mathbb{F}_p^\times$. This is enough to imply the congruence (12), in view of Stevens' calculation (via Dedekind sums) of the special values of the modular symbol associated to $\overline{\delta}_G$ (see [Ste82], Chapter 3).

In the following sections, we will fill in the details of the arguments outlined above. On the algebraic side, one very helpful fact is that the Selmer group $\text{Sel}_E(\mathbb{Q}_\infty)_p$ has a description which involves just the Galois module structure of $E[p^\infty]$. If $E$ has good, ordinary reduction at $p$, then $\text{Sel}_E(\mathbb{Q}_\infty)_p$ coincides with the Selmer group $\text{Sel}_{E[p^\infty]}(\mathbb{Q}_\infty)$ considered in [Gre89]. This simplifies the study of $\text{Sel}_E(\mathbb{Q}_\infty)[p]$ and is the basis for our proofs. Our results make sense in a rather general context – that of Selmer groups attached to modular forms. We will give the arguments in this more general setting. In particular, we will treat the case where $E$ has multiplicative reduction at $p$. If $E$ has nonsplit multiplicative reduction at $p$, then the arguments described in the introduction go through almost unchanged. However, if $E$ has split multiplicative at $p$, then there is an interesting difference: $\mathrm{Sel}_E(\mathbb{Q}_\infty)_p$ coincides with the "strict" Selmer group considered in [Gre89]. But $S_{E[p^\infty]}(\mathbb{Q}_\infty)$ is actually bigger. This corresponds to the fact that the associated $p$-adic L-function has a trivial zero. It is only when the trivial zero is included that our approach proceeds smoothly.

On the analytic side, the case of multiplicative reduction does not introduce any serious problems. The only difficulty appears in the comparison of the canonical periods and the Néron periods, where one has to generalize a theorem of Mazur (concerning the "Manin constant" of a strong Weil parametrization of an elliptic curve by $X_0(N)$) to the case of parametrizations by $X_1(N)$. This is mildly technical, as the reduction to characteristic $p$ of $X_1(N)$ can be quite complicated even when $N$ is divisible only by the first power of $p$, and Mazur's original argument does not generalize directly. Even though it is not strictly necessary, we have chosen parametrizations of elliptic curves by $X_1(N)$, rather than the more customary $X_0(N)$. The reason for this is connected to the behaviour of the $p$-adic L-functions under isogeny, and may be briefly described as follows. In the setting of Theorem (1.3), the elliptic curve corresponds to a nontrivial isogeny class $\mathcal{A}$. For the purposes of the proofs, one is forced to select a good representative of the class in question. The correct curve to use turns out to be the so-called *optimal* curve $E^{\mathrm{opt}}$ considered by Stevens in his paper [Ste89]. This curve is singled out by the requirement that there exist an embedding $E^{\mathrm{opt}} \hookrightarrow J_1(N)$, for the Jacobian variety $J_1(N)$ of $X_1(N)$. It seems that $E^{\mathrm{opt}}$ rather than any other (the strong Weil curve in $\mathcal{A}$ might be another candidate) is the "correct" curve to use when questions of integrality and congruence are as issue, as the associated lattice of Néron periods is conjectured to be "minimal" in a certain precise sense. Furthermore, parametrizations by $X_1(N)$ enjoy a certain universality property (see [Ste89], Theorem 1.9). Stevens has even conjectured that the optimal curve is characterized by *internal* rather than modular considerations. Namely, he conjectures that $E^{\mathrm{opt}}$ is the curve in $\mathcal{A}$ of minimal Faltings-Parshin height. We will in fact need some of the results of Stevens in our proofs of Lemma (3.6) and Corollary (3.8). However, we will not attempt to discuss the work of Stevens here. The reader will find a detailed account in the introduction of [Ste89]. (See also Remark (3.9) below.) Here we will content ourselves with a description of the situation when $\mathcal{A}$ is the isogeny class consisting of the three nonisomorphic curves of conductor 11. (This example was already mentioned above). In this case, the strong Weil curve is $X_0(11)$. As we have already remarked, the $\mu$-invariant of $X_0(11)$ is positive. The optimal curve is $X_1(11)$, and it turns out that the $\mu$-invariant vanishes. In view of Theorem 2.1 in [Ste85], one could even conjecture in general that the $\mu$-invariant vanishes for the optimal curve in an isogeny class.

Much of the research behind this article was carried out while the first author was supported by the National Science Foundation and by the American Institute of Mathematics. He is grateful to AIM for its support and hospitality during the Winter of 1998. The second author was supported by NSERC and also is grateful to AIM for the opportunity to visit in February, 1998.

## 2   Non-primitive Selmer groups

We will consider Selmer groups in a more general context than we actually need. Let $\Sigma$ be a finite set of primes of $\mathbb{Q}$ containing $p$ and $\infty$. Suppose that $\mathrm{Gal}\,(\mathbb{Q}_\Sigma/\mathbb{Q})$ acts continuously and linearly on a vector space $V_p$ over a field $F_p$. We assume that $d = \dim_{F_p}(V_p) < \infty$ and that $F_p$ is a finite extension of $\mathbb{Q}_p$. Let $\mathcal{O}$ denote the ring of integers of $F_p$. Let $T_p$ be a $\mathrm{Gal}\,(\mathbb{Q}_\Sigma/\mathbb{Q})$-invariant $\mathcal{O}$-lattice in $V_p$. Then $A = V_p/T_p$ is a discrete $\mathrm{Gal}\,(\mathbb{Q}_\Sigma/\mathbb{Q})$-module which is isomorphic to $(F_p/\mathcal{O})^d$ as an $\mathcal{O}$-module. If $d^{\pm}$ denotes the dimension of the $(\pm 1)$-eigenspaces for a complex conjugation, then $d = d^+ + d^-$. Fix an embedding $\overline{\mathbb{Q}} \to \overline{\mathbb{Q}}_p$. We can then identify $G_{\mathbb{Q}_p}$ with a decomposition group for some prime of $\overline{\mathbb{Q}}$ over $p$. We will assume that $V_p$ contains an $F_p$-subspace $W_p$ of dimension $d^+$ which is invariant under the action of $G_{\mathbb{Q}_p}$. Let $C$ denote the image of $W_p$ in $A$ under the canonical map $V_p \to A$, and let $D = A/C$.

The Selmer group $S_A(\mathbb{Q}_\infty)$ is defined by

$$S_A(\mathbb{Q}_\infty) = \ker\left(H^1(\mathbb{Q}_\Sigma/\mathbb{Q}_\infty, A) \to \prod_{\ell \in \Sigma} \mathcal{H}_\ell(\mathbb{Q}_\infty, A)\right)$$

where $\mathcal{H}_\ell(\mathbb{Q}_\infty, A)$ is defined as follows. If $\ell \neq p$, we simply let

$$\mathcal{H}_\ell(\mathbb{Q}_\infty, A) = \prod_{\eta | \ell} H^1((\mathbb{Q}_\infty)_\eta, A).$$

The product is over the finite set of primes $\eta$ of $\mathbb{Q}_\infty$ lying over $\ell$. There is unique prime $\eta_p$ of $\mathbb{Q}_\infty$ lying over $p$. Let $I_{\eta_p}$ denote the inertia subgroup of $G_{(\mathbb{Q}_\infty)_{\eta_p}}$. We define

$$\mathcal{H}_p(\mathbb{Q}_\infty, A) = H^1((\mathbb{Q}_\infty)_{\eta_p}, A)/L_{\eta_p}$$

where

$$L_{\eta_p} = \ker\left(H^1((\mathbb{Q}_\infty)_{\eta_p}, A) \to H^1(I_{\eta_p}, D)\right).$$

Thus if $\sigma$ is a 1-cocycle of $\mathrm{Gal}\,(\mathbb{Q}_\Sigma/\mathbb{Q}_\infty)$ with values in $A$, then its class $[\sigma]$ is in $S_A(\mathbb{Q}_\infty)$ if and only if $[\sigma|_{I_{\eta_p}}]$ is in the image of the map $H^1(I_{\eta_p}, C) \to H^1(I_{\eta_p}, A)$ and $[\sigma|_{G_{(\mathbb{Q}_\infty)_\eta}}] = 0$ for all $\eta|\ell, \ell \in \Sigma$ with $\ell \neq p$. Let $I_\eta$ denote the inertia subgroup of $G_{(\mathbb{Q}_\infty)_\eta}$. Then $G_{(\mathbb{Q}_\infty)_\eta}/I_\eta$ has profinite degree prime to $p$. So the last condition is equivalent to $[\sigma|_{I_\eta}] = 0$ for $\eta|\ell, \ell \in \Sigma$, and $\ell \neq p$. It is not hard to see that $S_A(\mathbb{Q}_\infty)$ is independent of the choice of $\Sigma$ (as long as $A$ is a $\mathrm{Gal}\,(\mathbb{Q}_\Sigma/\mathbb{Q})$-module and $p, \infty \in \Sigma$).

Now the groups $H^1(\mathbb{Q}_\Sigma/\mathbb{Q}_\infty, A)$, $H^2(\mathbb{Q}_\Sigma/\mathbb{Q}_\infty, A)$, $\mathcal{H}_\ell(\mathbb{Q}_\infty, A)$, and $S_A(\mathbb{Q}_\infty)$ are discrete $\mathcal{O}$-modules with a natural action of $\Gamma = \mathrm{Gal}\,(\mathbb{Q}_\infty/\mathbb{Q})$. Regarding them as $\Lambda$-modules, where $\Lambda = \mathcal{O}[[\Gamma]]$, they are known to be cofinitely generated. Using the results in sections 3 and 4 of [Gre89], one can easily verify the following statements:

1. $\mathrm{Corank}_\Lambda(H^1(\mathbb{Q}_\Sigma/\mathbb{Q}_\infty, A)) = d^- + \mathrm{Corank}_\Lambda(H^2(\mathbb{Q}_\Sigma/\mathbb{Q}_\infty, A))$.

2. $\mathrm{Corank}_\Lambda(\mathcal{H}_p(\mathbb{Q}_\infty)) = d^-$.

3. $\mathrm{Corank}_\Lambda(\mathcal{H}_\ell(\mathbb{Q}_\infty)) = 0$ if $\ell \neq p$.

In our results, we will generally assume that $S_A(\mathbb{Q}_\infty)$ is $\Lambda$-cotorsion. This clearly implies that $H^1(\mathbb{Q}_\Sigma/\mathbb{Q}_\infty, A)$ has $\Lambda$-corank $d^-$, that $H^2(\mathbb{Q}_\Sigma/\mathbb{Q}_\infty, A)$ is $\Lambda$-cotorsion, and that the cokernel of the map

$$\gamma : H^1(\mathbb{Q}_\Sigma/\mathbb{Q}_\infty, A) \to \prod_{\ell \in \Sigma} \mathcal{H}_\ell(\mathbb{Q}_\infty, A)$$

is $\Lambda$-cotorsion. Assuming that $p$ is odd, Proposition 4 of [Gre89] would then imply that $H^2(\mathbb{Q}_\Sigma/\mathbb{Q}_\infty, A) = 0$. (If $p = 2$ one would only get that $H^2(\mathbb{Q}_\Sigma/\mathbb{Q}, A)$ has exponent 2.) As for the cokernel of $\gamma$, the following result is crucial. Although proofs can be found elsewhere, we sketch a proof based on a generalization of a theorem of Cassels. We let $A^* = \mathrm{Hom}\,(T_p, \mu_{p^\infty})$, which is also a discrete $\mathcal{O}$-module with an action of $\mathrm{Gal}\,(\mathbb{Q}_\Sigma/\mathbb{Q})$.

**Proposition (2.1)** *Assume that $S_A(\mathbb{Q}_\infty)$ is $\Lambda$-cotorsion and that $H^0(\mathbb{Q}_\infty, A^*)$ is finite. Then $\gamma$ is surjective.*

*Proof.* It is enough to prove that $\mathrm{coker}(\gamma)$ is finite. The result would then follow because the Galois group $G_{(\mathbb{Q}_\infty)_\eta}$ has $p$-cohomological dimension 1 for any prime $\eta$ of $\mathbb{Q}_\infty$ and so $\mathcal{H}_\ell(\mathbb{Q}_\infty)$ is a divisible group for each $l$. Now let $\kappa : \Gamma \to 1 + p\mathbb{Z}_p$ be an isomorphism. For any $t \in \mathbb{Z}$, we let $A_t = A \otimes \kappa^t$, which is another $\mathrm{Gal}\,(\mathbb{Q}_\Sigma/\mathbb{Q})$-module. We can define a Selmer group $S_{A_t}(\mathbb{Q}_\infty)$ just as before. (For the local condition at $\eta_p$, one uses the $G_{\mathbb{Q}_p}$-invariant submodule $C_t = C \otimes \kappa^t$ of $A_t$ and the corresponding quotient $D_t = A_t/C_t$.) $S_{A_t}(\mathbb{Q}_\infty)$ is then the kernel of the map

$$\gamma_t : H^1(\mathbb{Q}_\Sigma/\mathbb{Q}_\infty, A_t) \to \prod_{\ell \in \Sigma} \mathcal{H}_\ell(\mathbb{Q}_\infty, A_t).$$

Clearly $A_t \cong A$ as $\mathrm{Gal}\,(\mathbb{Q}_\Sigma/\mathbb{Q}_\infty)$-modules and so it will suffice to show that $\mathrm{coker}(\gamma_t)$ is finite for at least one choice of $t$. We will do this by studying the cokernels of analogous "global-to-local" maps defined for the fields $\mathbb{Q}_n = \mathbb{Q}_\infty^{\Gamma_n}$ for all $n \geq 0$. Here $\Gamma_n = \Gamma^{p^n}$, $\mathbb{Q}_n$ is a cyclic extension of $\mathbb{Q}$ of degree $p^n$, and $\mathbb{Q}_\infty = \cup_n \mathbb{Q}_n$. Several requirements will be imposed on $t$ in the course of the proof. Two of the requirements are: (i) $(S_A(\mathbb{Q}_\infty) \otimes \kappa^t)^{\Gamma_n}$ is finite for all $n \geq 0$, (ii) $(A(\mathbb{Q}_\infty) \otimes \kappa^t)^{\Gamma_n}$ is finite for all $n \geq 0$. Since $S_A(\mathbb{Q}_\infty)$ is assumed to be $\Lambda$-cotorsion, it is easy to see that (i) is satisfied for all but finitely many values of $t$. The same is true for (ii) because $A(\mathbb{Q}_\infty) = H^0(\mathbb{Q}_\infty, A)$ is obviously $\Lambda$-cotorsion too.

Let $n$ and $t$ be fixed. We assume that $t$ satisfies the above requirements (i) and (ii). For brevity, we let $M = A_t$, $N = C_t$, and $K = \mathbb{Q}_n$. We define a Selmer group $S_M(K)$ by

$$S_M(K) = \ker\left(H^1(\mathbb{Q}_\Sigma/K, M) \to P_M(K)/L_M(K)\right).$$

Here $P_M(K) = \prod_\eta H^1(K_\eta, M)$, where $\eta$ runs over the primes of $K$ lying over those in $\Sigma$, and $L_M(K) = \prod_\eta L_\eta$ with $L_\eta = 0$ if $\eta \nmid p$ and $L_\eta = \ker\left(H^1(K_\eta, M) \to H^1(K_\eta, M/N)\right)$ for the unique prime $\eta$ of $K$ over $p$. It is easy to verify that the image of $S_M(K)$ under the restriction map $H^1(\mathbb{Q}_n, M) \to H^1(\mathbb{Q}_\infty, M)$ is contained in $S_{A_t}(\mathbb{Q}_\infty)^{\Gamma_n}$. We can identify $S_{A_t}(\mathbb{Q}_\infty)$ with $S_A(\mathbb{Q}_\infty) \otimes \kappa^t$ as $\Lambda$-modules, and so (i) implies that $S_{A_t}(\mathbb{Q}_\infty)^{\Gamma_n}$ is finite. The kernel of the restriction map is $H^1(\Gamma_n, M^{G_{\mathbb{Q}_\infty}})$. This has the same $\mathcal{O}$-corank as $H^0(\Gamma_n, M^{G_{\mathbb{Q}_\infty}}) = (A(\mathbb{Q}_\infty) \otimes \kappa^t)^{\Gamma_n}$, and so the kernel of the restriction map is finite by (ii). Thus $S_M(K)$ will be finite.

Now we will use the global duality theorems of Poitou and Tate. Let $U = \mathrm{Hom}\,(M, \mu_{p^\infty})$, which is a $\mathrm{Gal}\,(\mathbb{Q}_\Sigma/\mathbb{Q})$-module isomorphic to $\mathcal{O}^d$ as an $\mathcal{O}$-module. If $P_U(K) = \prod_\eta H^1(K_\eta, U)$, then local duality gives a perfect pairing

$$P_M(K) \times P_U(K) \to \mathbb{Q}_p/\mathbb{Z}_p. \tag{13}$$

We let $L_U(K)$ denote the orthogonal complement of $L_M(K)$ under the pairing (13). Then $L_U(K) = \prod_\eta L_\eta^\perp$, where $L_\eta^\perp$ is the orthogonal complement of $L_\eta$ under the local duality theorem for $K_\eta$. Note that if $\eta \nmid p$, then $L_\eta^\perp = H^1(K_\eta, U)$. Furthermore, let $G_M(K)$ and $G_U(K)$ denote the images of the maps

$$\alpha : H^1(\mathbb{Q}_\Sigma/K, M) \to P_M(K) \qquad \text{and} \qquad \beta : H^1(\mathbb{Q}_\Sigma/K, U) \to P_U(K)$$

respectively. Since $S_M(K) = \alpha^{-1}(L_M(K))$ is finite, it follows that $G_M(K) \cap L_M(K)$ is finite. Also, the global duality theorems assert that $G_M(K)$ and $G_U(K)$ are orthogonal complements under (13). Thus the cokernel of the map

$$H^1(\mathbb{Q}_\Sigma/K, M) \to P_M(K)/L_M(K) \tag{14}$$

is $P_M(K)/G_M(K)L_M(K)$ and this group is the Pontryagin dual of $G_U(K) \cap L_U(K)$.

The Euler characteristic for the $\mathrm{Gal}\,(\mathbb{Q}_\Sigma/K)$-module $M$ is

$$\sum_{i=0}^{2}(-1)^i \mathrm{corank}_\mathcal{O}(H^i(\mathbb{Q}_\Sigma/K, M)) = -d^- p^n$$

since $K = \mathbb{Q}_n$ is totally real and $[\mathbb{Q}_n : \mathbb{Q}] = p^n$. Hence the $\mathcal{O}$-corank of $H^1(\mathbb{Q}_\Sigma/K, M)$ is at least $d^- p^n$. On the other hand, if we exclude finitely many more values of $t$, then we can easily arrange for $P_M(K)/L_M(K)$ to have $\mathcal{O}$-corank equal to $\mathrm{corank}_\mathcal{O}(M/N)p^n = d^- p^n$. This is accomplished by making $H^0(K_\eta, M)$, $H^2(K_\eta, M)$ finite for all $\eta|l$, $l \in \Sigma$, $l \neq p$ and $H^0(K_\eta, M/N)$, $H^2(K_\eta, M/N)$ finite for the unique prime $\eta$ lying over $p$. Such a choice of $t$ is possible since $\kappa|_{G_{\mathbb{Q}_l}}$ has infinite order for each $\ell$ in the finite set $\Sigma$. Since $S_M(K)$ is finite, it now follows that the $\mathcal{O}$-corank of $H^1(\mathbb{Q}_\Sigma/K, M)$ is exactly $d^- p^n$ and that the cokernel of the map (14) is finite. Thus $G_U(K) \cap L_U(K)$ is a finite subgroup of $P_U(K)$. It also follows that $H^2(\mathbb{Q}_\Sigma/K, M)$ is finite (and even 0 if $p$ is odd, since $M$ is $\mathcal{O}$-divisible).

Let $S_U(K) = \beta^{-1}(L_U(K))$. Let $R_U^1(K) = \ker(\beta)$. By Poitou-Tate duality, $R_U^1(K)$ is the Pontryagin dual of $R_M^2(K) = \ker\left(H^2(\mathbb{Q}_\Sigma/K, M) \to \Pi_\eta H^2(K_\eta, M)\right)$. But this last group is clearly finite and so $R_U^1(K)$ is also finite (even 0 if $p$ is odd). Since $\beta(S_U(K)) = G_U(K) \cap L_U(K)$ is finite, it follows that $S_U(K)$ is a finite subgroup of $H^1(\mathbb{Q}_\Sigma/K, U)$. That is, $S_U(K) \subset H^1(\mathbb{Q}_\Sigma/K, U)_{\mathrm{tors}}$. We will show that

$$\#H^1(\mathbb{Q}_\Sigma/K, U)_{\mathrm{tors}} \leq \#H^0(\mathbb{Q}_\infty, A^*),$$

a bound independent of $n$. It will then follow that $\#(G_U(K) \cap L_U(K))$ is bounded as $n \to \infty$, and hence so is the cokernel of the map (14).

Let $V_p^* = \mathrm{Hom}\,(V_p, \mathbb{Q}_p(1))$ and $T_p^* = \mathrm{Hom}\,(T_p, \mathbb{Z}_p(1))$. Then one easily sees that $A^* = V_p^*/T_p^*$ and that $U = T_p^* \otimes \kappa^{-t}$. Put $M^* = A^* \otimes \kappa^{-t}$. Thus we have an exact sequence

$$0 \to U \to V_p^* \otimes \kappa^{-t} \to M^* \to 0.$$

Our hypothesis concerning $A^*$ implies that $H^0(\mathbb{Q}_\infty, M^*)$ is finite and $H^0(\mathbb{Q}_\infty, V_p^* \otimes \kappa^{-t}) = 0$. It follows that $H^1(\mathbb{Q}_\Sigma/K, U)_{\mathrm{tors}} \cong H^0(K, M^*)$, whose order is clearly bounded by the order of $H^0(\mathbb{Q}_\infty, A^*)$.

Thus, we have proved that for a suitable choice of $t$ (excluding just finitely many values of $t$), the cokernel of the map

$$\gamma'_{n,t} : H^1(\mathbb{Q}_\Sigma/\mathbb{Q}_n, A_t) \to P_{A_t}(\mathbb{Q}_n)/L_{A_t}(\mathbb{Q}_n)$$

is finite and of bounded order as $n \to \infty$. By taking the direct limit as $n \to \infty$, it follows that the cokernel of the limit map $\gamma'_t$ is finite. But since $A_t = A$ as a $\mathrm{Gal}\,(\mathbb{Q}_\Sigma/\mathbb{Q}_\infty)$-module

and $C_t = C$, it follows that the map

$$\gamma' : H^1(\mathbb{Q}_\Sigma/\mathbb{Q}_\infty, A) \to \left(H^1((\mathbb{Q}_\infty)_{\eta_p}, A/C) \times \prod_{\ell \in \Sigma, \ell \neq p} \mathcal{H}_\ell(\mathbb{Q}_\infty)\right)$$

has finite cokernel. In fact, the kernel of $\gamma'$ is the strict Selmer group $S_A^{\mathrm{str}}(\mathbb{Q}_\infty)$. It is clear that the cokernel of $\gamma$ is a homomorphic image of the cokernel of $\gamma'$. Hence the cokernel of $\gamma$ is finite, and, as remarked earlier, $\gamma$ must be surjective. ∎

**Remark (2.2)** It will be useful later to point out one other consequence of the above proof. We make the same requirements on the choice of $t$. In particular, we have that $H^0(K_\eta, M)$, $H^2(K_\eta, M)$ are both finite for all $\eta | l$, $l \in \Sigma$, $l \neq p$. Then $H^1(K_\eta, M)$ is finite too. Assuming that $\Sigma - \{p, \infty\}$ is nonempty, choose one such prime $l_0$. Suppose that in defining $S_M(K)$ we take $L_\eta = H^1(K_\eta, M)$ for all $\eta | l_0$ (instead of taking $L_\eta = 0$). Then $L_\eta^\perp = 0$. Assume that $t$ is also chosen so that $H^0(K_\eta, M^*)$ is finite. (This of course implies that $H^0(K, M^*)$ is finite.) With these changes, it then turns out that $S_U(K) = 0$ and hence the map (14) is now surjective. This becomes clear from the following commutative diagram:

$$\begin{array}{ccc} H^0(K, M^*) & \xrightarrow{\sim} & H^1(\mathbb{Q}_\Sigma/K, U)_{\mathrm{tors}} \\ \downarrow & & \downarrow \\ H^0(K_\eta, M^*) & \xrightarrow{\sim} & H^1(K_\eta, U)_{\mathrm{tors}} \end{array}$$

Here $\eta$ is any prime of $K$ dividing $l_0$. The horizontal isomorphisms are coboundary maps. The first vertical map is clearly injective, and therefore so is the second. This shows that $S_U(K) \subset H^1(\mathbb{Q}_\Sigma/K, U)_{\mathrm{tors}}$ must indeed be zero.

Let $\Sigma_0$ be any finite set of primes of $\mathbb{Q}$ which does not contain $p$ or $\infty$. Choose the set $\Sigma$ large enough so that $\Sigma_0 \subset \Sigma$. The non-primitive Selmer group for $A$ and $\Sigma_0$ is defined by

$$S_A^{\Sigma_0}(\mathbb{Q}_\infty) = \ker\left(H^1(\mathbb{Q}_\Sigma/\mathbb{Q}_\infty, A) \to \prod_{\ell \in \Sigma - \Sigma_0} \mathcal{H}_\ell(\mathbb{Q}_\infty)\right).$$

Obviously, $S_A(\mathbb{Q}_\infty) \subset S_A^{\Sigma_0}(\mathbb{Q}_\infty)$. Proposition (2.1), together with the fact that $\mathcal{H}_\ell(\mathbb{Q}_\infty)$ is $\Lambda$-cotorsion for $\ell \neq p$, immediately gives the following result:

**Corollary (2.3)** *With assumptions as in (2.1), we have*

$$S_A^{\Sigma_0}(\mathbb{Q}_\infty)/S_A(\mathbb{Q}_\infty) \cong \prod_{\ell \in \Sigma_0} \mathcal{H}_\ell(\mathbb{Q}_\infty)$$

as $\Lambda$-modules. Furthermore, $S_A^{\Sigma_0}(\mathbb{Q}_\infty)$ is $\Lambda$-cotorsion and

$$\mathrm{corank}_\mathcal{O}(S_A^{\Sigma_0}(\mathbb{Q}_\infty)) = \mathrm{corank}_\mathcal{O}(S_A(\mathbb{Q}_\infty)) + \sum_{\ell \in \Sigma_0} \mathrm{corank}_\mathcal{O}(\mathcal{H}_\ell(\mathbb{Q}_\infty)).$$

The $\mu$-invariants of $S_A^{\Sigma_0}(\mathbb{Q}_\infty)\hat{\phantom{o}}$ and $S_A(\mathbb{Q}_\infty)\hat{\phantom{o}}$ are equal.

The structure of $\mathcal{H}_\ell(\mathbb{Q}_\infty)$ can be studied using Proposition 2 of [Gre89]. Since we will be primarily interested in its $\mathcal{O}$-corank, we discuss that first. Let $s_\ell$ denote the number of primes $\eta$ of $\mathbb{Q}_\infty$ lying over $\ell$. That is, $s_\ell = [\Gamma : \Gamma_\ell]$, where $\Gamma_\ell$ denotes the decomposition subgroup of $\Gamma$ for any such $\eta$. If $\gamma_\ell$ denotes the corresponding Frobenius automorphism in $\Gamma$, then $\Gamma_\ell$ is generated topologically by $\gamma_\ell$. One can determine $s_\ell$ quite easily. For an odd prime $p$, $s_\ell$ is the largest power of $p$ such that $\ell^{p-1} \equiv 1 \pmod{ps_\ell}$. Now let $d_\ell = \dim_{F_p}(H^0((\mathbb{Q}_\infty)_\eta, V_p^*))$, where $V_p^* = \mathrm{Hom}\,(V_p, \mathbb{Q}_p(1))$ as before. Proposition 2 of [Gre89] easily implies that $\mathrm{corank}_\mathcal{O}(\mathcal{H}_\ell(\mathbb{Q}_\infty)) = s_\ell d_\ell$.

The value of $d_\ell$ is also not hard to determine. Let $\mathrm{Frob}_\ell$ denote the Frobenius automorphism in $\mathrm{Gal}\,(\mathbb{Q}_\ell^{\mathrm{unr}}/\mathbb{Q}_\ell)$. (Thus, $\gamma_\ell$ is the restriction of $\mathrm{Frob}_\ell$ to $(\mathbb{Q}_\infty)_\eta$.) Let $\alpha_1, \alpha_2, \ldots, \alpha_{e_\ell}$ denote the eigenvalues of $\mathrm{Frob}_\ell$ (counting multiplicities) acting on $(V_p)_{I_\ell}$. Here $I_\ell$ is the inertia subgroup $\mathrm{Gal}\,(\overline{\mathbb{Q}}_\ell/\mathbb{Q}_\ell^{\mathrm{unr}})$ of $G_{\mathbb{Q}_\ell}$, $(V_p)_{I_\ell}$ denotes the maximal quotient of $V_p$ on which $I_\ell$ acts trivially, and $e_\ell = \dim_{F_p}((V_p)_{I_\ell})$. Since $(\mathbb{Q}_\infty)_\eta \subset \mathbb{Q}_\ell^{\mathrm{unr}}$ and $I_\ell$ acts trivially on $\mathbb{Q}_p(1)$, we have

$$H^0((\mathbb{Q}_\infty)_\eta, V_p^*) \subset H^0(I_\ell, V_p^*) = \mathrm{Hom}\,_{F_p}((V_p)_{I_\ell}, \mathbb{Q}_p(1)).$$

The eigenvalues of $\mathrm{Frob}_\ell$ acting on this last vector space are $\ell\alpha_1^{-1}, \ldots, \ell\alpha_{e_\ell}^{-1}$. Noting that the action of $\mathrm{Gal}\,(\mathbb{Q}_\ell^{\mathrm{unr}}/(\mathbb{Q}_\infty)_\eta)$ must be through a finite group of order prime to $p$, one sees easily that $d_\ell$ is the number of $i$'s such that $\ell\alpha_i^{-1}$ is a principal unit in $F_p(\alpha_i)$. These values of $\ell\alpha_i^{-1}$ are precisely the eigenvalues of $\gamma_\ell$ acting on $H^0((\mathbb{Q}_\infty)_\eta, V_p^*)$, counting multiplicities. Alternatively, we can describe $d_\ell$ in terms of the element

$$\mathscr{P}_\ell = \prod_{i=1}^{e_\ell}(1 - \alpha_i\ell^{-1}\gamma_\ell) \in \mathcal{O}[[\Gamma_\ell]]. \tag{15}$$

Identifying $\mathcal{O}[[\Gamma_\ell]]$ with the power series ring $\mathcal{O}[[T_\ell]]$, where $T_\ell = \gamma_\ell - 1$, we can then factor $\mathscr{P}_\ell = \mathscr{P}_\ell(T_\ell)$ as a product of a power of $\pi$ (a uniformizing parameter for $\mathcal{O}$), an invertible power series, and a distinguished polynomial. The power of the uniformizing parameter is 1. (That is, the $\mu$-invariant is 0.) The degree of the distinguished polynomial is $d_\ell$. If we view $\mathscr{P}_\ell$ as an element of $\mathcal{O}[[\Gamma]] = \mathcal{O}[[T]]$, with $T = \gamma - 1$ for a fixed topological generator of $\Gamma$, then $\mathscr{P}_l$ will still have $\mu$-invariant 0 and its distinguished polynomial factor, which we denote by $h_\ell(T)$, will have degree $s_\ell d_\ell$. Thus the $\mathcal{O}$-corank of $\mathcal{H}_\ell(\mathbb{Q}_\infty)$ is $\deg(h_\ell(T))$.

The polynomial $h_\ell(T)$ generates the characteristic ideal of the $\Lambda$-module $\mathcal{H}_\ell(\mathbb{Q}_\infty)\hat{\,}$. To see this, it is enough to verify that $\mathscr{P}_\ell$ generates the characteristic ideal of the $\mathcal{O}[[\Gamma_\ell]]$-module $H^1((\mathbb{Q}_\infty)_\eta, A)\hat{\,}$, as $\mathcal{H}_\ell(\mathbb{Q}_\infty)\hat{\,}$ is obtained by tensoring with $\Lambda$ over $\mathcal{O}[[\Gamma_\ell]]$. The inertia subgroup $I_\ell$ of $G_{\mathbb{Q}_\ell}$ contains a unique subgroup $J_\ell$ such that $I_\ell/J_\ell \cong \mathbb{Z}_p$. Also, $J_\ell$ has profinite order prime to $p$ and $I_\ell/J_\ell \cong \mathbb{Z}_p(1)$ for the natural action of $\text{Gal}\,(\mathbb{Q}_\ell^{\text{unr}}/\mathbb{Q}_\ell)$. Let $G = \text{Gal}\,(\overline{\mathbb{Q}}_\ell^{J_\ell}/(\mathbb{Q}_\infty)_\eta) = G_{(\mathbb{Q}_\infty)_\eta}/I_\ell$. Since $H^1(J_\ell, A) = 0$, we have $H^1((\mathbb{Q}_\infty)_\eta, A) = H^1(G, A^{J_\ell})$. Also, $A^{J_\ell} \cong A_{J_\ell}$ canonically. If we let $\overline{I}_\ell = I_\ell/J_\ell$, then $G/\overline{I}_\ell$ has profinite order prime to $p$ and so $H^1(G, A_{I_\ell}) \cong H^1(\overline{I}_\ell, A_{J_\ell})^{G/\overline{I}_\ell}$. If we let $\epsilon_\ell$ be a topological generator of $\overline{I}_\ell = \mathbb{Z}_p(1)$, then

$$H^1(\overline{I}_l, A_{J_\ell}) = H^1(\overline{I}_\ell, A_{J_\ell}/(\epsilon_\ell - 1)A_{J_\ell}) \cong \text{Hom}\,(\mathbb{Z}_\ell(1), A_{I_\ell}) \cong A_{I_\ell}(-1),$$

where the isomorphisms are equivariant for the action of $\text{Gal}\,(\mathbb{Q}_\ell^{\text{unr}}/\mathbb{Q}_\ell)$. The eigenvalues of $\text{Frob}_\ell$ acting on $A_{I_\ell}(-1)\hat{\,}$ are the numbers $\ell\alpha_i^{-1}, 1 \leq i \leq e_\ell$, and those numbers which are principal units are the eigenvalues of $\gamma_\ell$ acting on $(A_{I_\ell}(-1)^{G/\overline{I}_\ell})\hat{\,}$, again using the fact that $G/\overline{I}_\ell$ has profinite order prime to $p$. These remarks imply that $\mathscr{P}_\ell$ does generate the characteristic ideal of $H^1((\mathbb{Q}_\infty)_\eta, A)\hat{\,}$ as an $\mathcal{O}[[\Gamma_\ell]]$-module. We have proved the following results:

**Proposition (2.4)** *Let $P_\ell(X) = \det((1 - \text{Frob}_\ell X)|_{(V_p)_{I_\ell}}) \in \mathcal{O}[X]$. Let $\mathscr{P}_\ell = P_\ell(\ell^{-1}\gamma_\ell) \in \Lambda = \mathcal{O}[[\Gamma]]$, where $\gamma_\ell$ denotes the Frobenius automorphism for $\ell$ in $\Gamma = \text{Gal}\,(\mathbb{Q}_\infty/\mathbb{Q})$. The characteristic ideal of the $\Lambda$-module $\mathcal{H}_\ell(\mathbb{Q}_\infty)\hat{\,}$ is generated by $\mathscr{P}_\ell$. Its $\mu$-invariant is zero. Its $\lambda$-invariant is equal to $s_\ell d_\ell$. Here $s_\ell$ is the largest power of $p$ dividing $(\ell^{p-1} - 1)/p$ and $d_\ell$ is the multiplicity of $X = \widetilde{\ell}^{-1}$ as a root of $\widetilde{P}_\ell(X) \in k[X]$, where $k$ is the residue field of $\mathcal{O}$, and the $\sim$ means reduction modulo $\mathfrak{m}$, where $\mathfrak{m} = (\pi)$ denotes the maximal ideal of $\mathcal{O}$.*

Remark. $\mathscr{P}_l$ satisfies an interpolation property involving Euler factors evaluated at $s = 1$ (if $V_p$ arises from a compatible system of $l$-adic (or $\lambda$-adic) representations of $G_\mathbb{Q}$). Namely, let $\rho$ be any character of $\Gamma$. Then $\rho(\mathscr{P}_l) = P_l(\rho(l)l^{-1})$, regarding $\rho$ as a Dirichlet character.

Now we will discuss the nonexistence of nonzero finite $\Lambda$-submodules in the Pontrjagin duals of Selmer groups. One finds results about this in [Gre99] for the case of $\text{Sel}_E(F_\infty)_p$, where $E$ is an elliptic curve defined over a number field $F$ and $F_\infty$ is the cyclotomic $\mathbb{Z}_p$-extension of $F$. (See Propositions 4.14 and 4.15 in [Gre99].) Here we will prove a much easier result, which will be sufficient for our purposes. Recall that $D = A/C$ is an $\mathcal{O}$-module on which $G_{\mathbb{Q}_p}$ acts. We let $\text{Ram}(A)$ denote the set of primes $\ell$ ($\neq p, \infty$) such that the action of $G_{\mathbb{Q}_\ell}$ on $A$ is ramified.

**Proposition (2.5)** *Let $p$ be an odd prime. Assume that $S_A(\mathbb{Q}_\infty)$ is $\Lambda$-cotorsion and that $D$ is unramified for the action of $G_{\mathbb{Q}_p}$. Suppose that $\Sigma_0$ is a subset of $\Sigma - \{p, \infty\}$ which contains $\mathrm{Ram}(A)$. Then $S_A^{\Sigma_0}(\mathbb{Q}_\infty)\hat{\ }$ has no nonzero, finite $\Lambda$-submodules.*

*Proof.* $S_A^{\Sigma_0}(\mathbb{Q}_\infty)$ doesn't depend on the choice of $\Sigma$, as long as $\Sigma_0 \cup \{p, \infty\}$ is contained in $\Sigma$. In this proof, we may therefore take $\Sigma = \Sigma_0 \cup \{p, \infty\}$. Then, by definition, $S_A^{\Sigma_0}(\mathbb{Q}_\infty) = \ker(H^1(\mathbb{Q}_\Sigma/\mathbb{Q}_\infty, A) \to \mathcal{H}_p(\mathbb{Q}_\infty))$. Since $S_A(\mathbb{Q}_\infty)$ is $\Lambda$-cotorsion, so is $S_A^{\Sigma_0}(\mathbb{Q}_\infty)$. Both $H^1(\mathbb{Q}_\Sigma/\mathbb{Q}_\infty, A)$ and $\mathcal{H}_p(\mathbb{Q}_\infty)$ have $\Lambda$-corank equal to $d^-$. We will show that $\mathcal{H}_p(\mathbb{Q}_\infty)$ is $\Lambda$-cofree. It is then clear that the map defining $S_A^{\Sigma_0}(\mathbb{Q}_\infty)$ is surjective and so

$$H^1(\mathbb{Q}_\Sigma/\mathbb{Q}_\infty, A)/S_A^{\Sigma_0}(\mathbb{Q}_\infty) \cong \mathcal{H}_p(\mathbb{Q}_\infty)$$

as $\Lambda$-modules. (This would usually follow from Proposition 2.1, even without knowing something about $\mathcal{H}_p(\mathbb{Q}_\infty)$.) Proposition 5 of [Gre89] asserts that $H^1(\mathbb{Q}_\Sigma/\mathbb{Q}_\infty, A)\hat{\ }$ has no nonzero, finite $\Lambda$-submodules because we know that $H^2(\mathbb{Q}_\Sigma/\mathbb{Q}_\infty, A) = 0$. (We are assuming that $p$ is odd.) The assertion about $(S_A^{\Sigma_0}(\mathbb{Q}_\infty))\hat{\ }$ is then a consequence of the following lemma (whose simple proof can be found on page 123 of [Gre89]:

**Lemma (2.6)** *Let $Y$ be a finitely generated $\Lambda$-module, $Z$ a free $\Lambda$-submodule. If $Y$ contains no nonzero, finite $\Lambda$-submodule, then the same is true for $Y/Z$.*

We just apply this to $Y = H^1(\mathbb{Q}_\Sigma/\mathbb{Q}_\infty, A)\hat{\ }, Z = \mathcal{H}_p(\mathbb{Q}_\infty)\hat{\ }$, and $Y/Z = S_A^{\Sigma_0}(\mathbb{Q}_\infty)\hat{\ }$.

It remains to show that $\mathcal{H}_p(\mathbb{Q}_\infty)$ is $\Lambda$-cofree when $D$ is unramified. We first verify that $H^1(\mathbb{Q}_p, D)$ is $\mathcal{O}$-cofree. The exact sequence $0 \to D[\pi] \to D \xrightarrow{\pi} D \to 0$ induces an injective map $H^1(\mathbb{Q}_p, D)/\pi H^1(\mathbb{Q}_p, D) \to H^2(\mathbb{Q}_p, D[\pi])$. This last group is dual to $H^0(\mathbb{Q}_p, \mathrm{Hom}(D[\pi], \mu_p))$ which is obviously trivial since $p$ is an odd prime. Thus $H^1(\mathbb{Q}_p, D)$ is a divisible $\mathcal{O}$-module. Its Pontryagin dual is a torsion-free, finitely generated $\mathcal{O}$-module and must therefore be free. Thus $H^1(\mathbb{Q}_p, D)$ is $\mathcal{O}$-cofree. Its $\mathcal{O}$-corank is equal to $d^- + \mathrm{corank}_\mathcal{O}(H^0(\mathbb{Q}_p, D))$ since $D$ has $\mathcal{O}$-corank $d^-$ and $H^2(\mathbb{Q}_p, D) = 0$.

Now by the inflation-restriction exact sequence together with the fact that $\Gamma$ has $p$-cohomological dimension 1, the restriction map $H^1(\mathbb{Q}_p, D) \to H^1((\mathbb{Q}_\infty)_{\eta_p}, D)^\Gamma$ is surjective, where $\Gamma$ is identified with $\mathrm{Gal}((\mathbb{Q}_\infty)_{\eta_p}/\mathbb{Q}_p)$. The kernel of this map is $H^1(\Gamma, D^{G_{(\mathbb{Q}_\infty)_{\eta_p}}})$, which is easily seen to have $\mathcal{O}$-corank equal to that of $H^0(\mathbb{Q}_p, D)$. Thus $H^1((\mathbb{Q}_\infty)_{\eta_p}, D)^\Gamma$ is $\mathcal{O}$-cofree and has $\mathcal{O}$-corank $d^-$. By proposition 1 of [Gre89], we see that $H^1((\mathbb{Q}_\infty)_{\eta_p}, D)$ has $\Lambda$-corank $d^-$. Hence $X = H^1((\mathbb{Q}_\infty)_{\eta_p}, D)\hat{\ }$ is a $\Lambda$-module of rank $d^-$ such that $X/TX$ is $\Lambda/T\Lambda$-free of rank $d^-$. A simple use of Nakayama's lemma shows that $X$ is a free $\Lambda$-module.

It follows from the definition of $\mathcal{H}_p(\mathbb{Q}_\infty)$ that we have an isomorphism

$$\mathcal{H}_p(\mathbb{Q}_\infty) \cong \mathrm{im}\,(H^1((\mathbb{Q}_\infty)_{\eta_p}, A) \to H^1(I_{\eta_p}, D))$$

as $\Lambda$-modules. Since $G_{(\mathbb{Q}_\infty)_{\eta_p}}$ has $p$-cohomological dimension equal to 1, we also have $H^2((\mathbb{Q}_\infty)_{\eta_p}, C) = 0$. Hence the exact sequence $0 \to C \to A \to D \to 0$ implies that the map $H^1((\mathbb{Q}_\infty)_{\eta_p}, A) \to H^1((\mathbb{Q}_\infty)_{\eta_p}, D)$ is surjective. Also, $\ker(H^1((\mathbb{Q}_\infty)_{\eta_p}, D) \to H^1(I_{\eta_p}, D))$ is isomorphic to $H^1(G_{(\mathbb{Q}_\infty)_{\eta_p}}/I_{\eta_p}, D)$ because $I_{\eta_p}$ acts trivially on $D$. Since $G_{(\mathbb{Q}_\infty)_{\eta_p}}/I_{\eta_p}$ is topologically cyclic, this kernel is a quotient of $D$ and so is $\mathcal{O}$-cofree. (It would be sufficient to assume just that $D^{I_{\eta_p}}$ is $\mathcal{O}$-cofree.) Thus $\mathcal{H}_p(\mathbb{Q}_\infty)\hat{\,}$ is isomorphic to a $\Lambda$-submodule $Y$ of $X = H^1((\mathbb{Q}_\infty)_{\eta_p}, D)\hat{\,}$ such that $X/Y$ is $\mathcal{O}$-cofree. Since $X$ is $\Lambda$-free, one sees that $Y$ must be reflexive and hence also free as a $\Lambda$-module. Hence $\mathcal{H}_p(\mathbb{Q}_\infty)$ is indeed $\Lambda$-cofree. ∎

**Remark (2.7)** In proposition 2.5, it is not necessary to assume that $D$ is unramified for the action of $G_{\mathbb{Q}_p}$. In the above proof, only two properties of $D$ were actually used. First, that $\text{Hom}_{G_{\mathbb{Q}_p}}(D[\pi], \mu_p) = 0$, which suffices to show that $H^1(\mathbb{Q}_p, D)$ is $\mathcal{O}$-cofree and that $H^2(\mathbb{Q}_p, D) = 0$. Second, that $D^{I_{\eta_p}}$ is $\mathcal{O}$-cofree which is used at the end of the proof. This result covers virtually all the cases we are interested in. However, we will sketch another approach which gives the same conclusion with slightly different hypotheses. Assume that $A$ satisfies the hypotheses of proposition 2.1, that $\Sigma_0$ is nonempty, and that $H^1(G_{(\mathbb{Q}_\infty)_{\eta_p}}/I_{\eta_p}, D^{I_{\eta_p}}) = 0$. (Note that the last hypothesis is valid if $D^{I_{\eta_p}} = 0$.) Then, according to remark 2.2, the following map can be assumed to be surjective:

$$H^1(\mathbb{Q}_\Sigma/\mathbb{Q}, M) \to H^1(\mathbb{Q}_p, M/N)$$

where $M = A \otimes \kappa^t$ for a suitable $t$. (We are assuming that in addition to other requirements, $t$ is also chosen so that $H^2(\mathbb{Q}_p, N) = 0$, so that the map $H^1(\mathbb{Q}_p, M) \to H^1(\mathbb{Q}_p, M/N)$ is surjective. Here $N = C \otimes \kappa^t$ and such a choice is possible.) Using the inflation-restriction sequence and the fact that $\Gamma$ has cohomological dimension 1, one sees that the map

$$H^1(\mathbb{Q}_\Sigma/\mathbb{Q}_\infty, M)^\Gamma \to H^1((\mathbb{Q}_\infty)_{\eta_p}, M/N)^\Gamma$$

is also surjective. Now $S_M^{\Sigma_0}(\mathbb{Q}_\infty) \cong S_A^{\Sigma_0}(\mathbb{Q}_\infty) \otimes \kappa^t$ as $\Lambda$-modules and the exact sequence

$$0 \to S_M^{\Sigma_0}(\mathbb{Q}_\infty) \to H^1(\mathbb{Q}_\Sigma/\mathbb{Q}_\infty, M) \to H^1((\mathbb{Q}_\infty)_{\eta_p}, M/N) \to 0$$

together with the snake lemma imply that the map $S_M^{\Sigma_0}(\mathbb{Q}_\infty)_\Gamma \to H^1(\mathbb{Q}_\Sigma/\mathbb{Q}_\infty, M)_\Gamma$ is injective. Since $H^1(\mathbb{Q}_\Sigma/\mathbb{Q}_\infty, M)\hat{\,}$ has no nonzero finite $\Lambda$-submodule, it follows that, if $t$ is again chosen suitably, $H^1(\mathbb{Q}_\Sigma/\mathbb{Q}_\infty, M)_\Gamma = 0$. Therefore, $S_M^{\Sigma_0}(\mathbb{Q}_\infty)_\Gamma = 0$ and the conclusion in proposition 2.5 follows immediately.

The next proposition allows one to determine the $\lambda$-invariant of $S_A^{\Sigma_0}(\mathbb{Q}_\infty)$ in terms of the Galois module $A[\pi]$, under certain hypotheses. Recall that $\pi$ is a generator of the maximal ideal of $\mathcal{O}$. We define a Selmer group for $A[\pi]$ in the following way. Consider the exact sequence $0 \to C[\pi] \to A[\pi] \to D[\pi] \to 0$ of $G_{\mathbb{Q}_p}$-modules. For any subset $\Sigma_0$ of $\Sigma - \{p, \infty\}$, let

$$S_{A[\pi]}^{\Sigma_0}(\mathbb{Q}_\infty) = \ker\left(H^1(\mathbb{Q}_\Sigma/\mathbb{Q}_\infty, A[\pi]) \to \prod_{\ell \in \Sigma - \Sigma_0} \mathcal{H}_\ell(\mathbb{Q}_\infty, A[\pi])\right).$$

where for $\ell \neq p$, we define $\mathcal{H}_\ell(\mathbb{Q}_\infty, A[\pi]) = \prod_{\eta | \ell} H^1(I_\eta, A[\pi])$ and, for $\ell = p$, we define $\mathcal{H}_p(\mathbb{Q}_\infty, A[\pi]) = H^1(I_{\eta_p}, D[\pi])$. With this definition, which is entirely analogous to the definition of $S_A^{\Sigma_0}(\mathbb{Q}_\infty)$, we can prove the following useful result.

**Proposition (2.8)** *Let $p$ be an odd prime. Assume that $\Sigma_0$ is a subset of $\Sigma - \{p, \infty\}$ containing $\mathrm{Ram}(A)$. Assume that $I_{\eta_p}$ acts trivially on $D$ and that $H^0(\mathbb{Q}, A[\pi]) = 0$. Then*

$$S_A^{\Sigma_0}(\mathbb{Q}_\infty)[\pi] \cong S_{A[\pi]}^{\Sigma_0}(\mathbb{Q}_\infty).$$

*Consequently, $S_A(\mathbb{Q}_\infty)$ is $\Lambda$-cotorsion and has $\mu$-invariant zero if and only if $S_{A[\pi]}^{\Sigma_0}(\mathbb{Q}_\infty)$ is finite. If this is so, then the $\lambda$-invariant of $S_A^{\Sigma_0}(\mathbb{Q}_\infty)$ is equal to $\dim_{\mathcal{O}/\mathfrak{m}}(S_{A[\pi]}^{\Sigma_0}(\mathbb{Q}_\infty))$.*

*Proof.* Since $H^0(\mathbb{Q}, A[\pi]) = 0$ and $\Gamma$ is a pro-$p$ group, it follows that $H^0(\mathbb{Q}_\infty, A) = 0$. The natural map

$$H^1(\mathbb{Q}_\Sigma/\mathbb{Q}_\infty, A[\pi]) \to H^1(\mathbb{Q}_\Sigma/\mathbb{Q}_\infty, A)[\pi]$$

induced from the exact sequence $0 \to A[\pi] \to A \xrightarrow{\pi} A \to 0$ is therefore an isomorphism. We must compare the local conditions defining $S_A^{\Sigma_0}(\mathbb{Q}_\infty)[\pi]$ and $S_{A[\pi]}^{\Sigma_0}(\mathbb{Q}_\infty)$. Suppose that $\sigma$ is a 1-cocycle of $\mathrm{Gal}(\mathbb{Q}_\Sigma/\mathbb{Q}_\infty)$ with values in $A[\pi]$. First consider the local condition at $\eta | \ell$, where $\ell \neq p$ and $\ell \in \Sigma - \Sigma_0$. Then $I_\eta$ acts trivially on $A$. The map $H^1(I_\eta, A[\pi]) \to H^1(I_\eta, A)$ is injective because $H^0(I_\eta, A) = A$ is divisible. Thus, the local conditions at $\ell$ defining $S_A^{\Sigma_0}(\mathbb{Q}_\infty)[\pi]$ and $S_{A[\pi]}^{\Sigma_0}(\mathbb{Q}_\infty)$ are equivalent. Also, the map

$$H^1(I_{\eta_p}, D[\pi]) \to H^1(I_{\eta_p}, D)$$

is injective because $H^0(I_{\eta_p}, D) = D$ is divisible. Hence the local conditions on $\sigma$ defining the two Selmer groups are equivalent, and this proves the equality.

Now $S_A(\mathbb{Q}_\infty)$ and $S_A^{\Sigma_0}(\mathbb{Q}_\infty)$ have the same $\Lambda$-corank and the same $\mu$-invariant if they are $\Lambda$-cotorsion. Obviously, $S_A^{\Sigma_0}(\mathbb{Q}_\infty)$ is $\Lambda$-cotorsion and has $\mu$-invariant 0 if and only if $S_A^{\Sigma_0}(\mathbb{Q}_\infty)[\pi]$ is finite, which in turn is equivalent to the finiteness of $S_{A[\pi]}^{\Sigma_0}(\mathbb{Q}_\infty)$. Assuming this is so, $S_A^{\Sigma_0}(\mathbb{Q}_\infty)\hat{\,}$ would be a finitely generated $\mathcal{O}$-module. By Proposition (2.5) its $\mathcal{O}$-torsion submodule is 0, and so $S_A^{\Sigma_0}(\mathbb{Q}_\infty)$ is $\mathcal{O}$-divisible. That is, $S_A^{\Sigma_0}(\mathbb{Q}_\infty) \cong (F_p/\mathcal{O})^\lambda$, where $\lambda = \mathrm{corank}_\mathcal{O}(S_A^{\Sigma_0}(\mathbb{Q}_\infty)) \geq 0$. It is clear that $\lambda = \dim_{\mathcal{O}/\mathfrak{m}}(S_A^{\Sigma_0}(\mathbb{Q}_\infty)) = \dim_{\mathcal{O}/\mathfrak{m}} S_{A[\pi]}^{\Sigma_0}(\mathbb{Q}_\infty)$, as stated. ∎

**Remark (2.9)** The hypothesis that $I_{\eta_p}$ acts trivially on $D$ is adequate for the theorems stated in the introduction. But all that is needed in the proof is that $H^0(I_{\eta_p}, D)$ is divisible. Thus, if $D^{I_{\eta_p}} = 0$, the conclusion is still true.

At this point we can justify most of the steps in the proof of theorem (1.4) outlined in the introduction. If $E$ is a modular elliptic curve defined over $\mathbb{Q}$, we take $T_p = T_p(E)$, the Tate module for $E$, and $V_p = T_p(E) \otimes_{\mathbb{Z}_p} \mathbb{Q}_p$, a 2-dimensional $\mathbb{Q}_p$-representation space for $G_{\mathbb{Q}}$. Then $A = V_p/T_p$ is isomorphic to $E[p^\infty]$ as a $G_{\mathbb{Q}}$-module. The ring $\mathcal{O}$ is $\mathbb{Z}_p$, and $d^+ = d^- = 1$. Assume that $E$ has good ordinary reduction at $p$. Then viewing $A$ as a $G_{\mathbb{Q}_p}$-module, we define $C = \ker(E[p^\infty] \to \widetilde{E}[p^\infty])$, where $\widetilde{E}$ is the reduction of $E$ modulo $p$. Then $D \cong \widetilde{E}[p^\infty]$ is unramified as a $G_{\mathbb{Q}_p}$-module and is isomorphic to $\mathbb{Q}_p/\mathbb{Z}_p$ as a $\mathbb{Z}_p$-module. By the Weil pairing, one sees that the inertia group $I_p$ acts on $C$ by the $p$-cyclotomic character $\chi$. That is, $C \cong \mu_{p^\infty}$ as an $I_p$-module. Note that the action of $G_{\mathbb{Q}_p}$ on $A$ depends on fixing an embedding of the field $\mathbb{Q}(E[p^\infty])$ into $\overline{\mathbb{Q}_p}$. The subgroup $C$ of $A$ depends on the corresponding choice of a prime of $\mathbb{Q}(E[p^\infty])$ lying above $p$. Having fixed such a prime, the subgroup $C$ is determined by the action of $I_p$. Also, since we are assuming that $p$ is odd, the subgroup $C[p]$ of $A[p]$ is determined by the action of $I_p : C[p] \cong \mu_p$ and $D[p]$ is the maximal quotient of $A[p]$ on which $I_p$ acts trivially.

In [Gre99], one can find a proof that $\operatorname{im}(\kappa_{\eta_p}) = L_{\eta_p}$. This result, together with the fact that $\operatorname{im}(\kappa_\eta) = 0$ for all primes $\eta$ of $\mathbb{Q}_\infty$ not lying over $p$, implies that $\operatorname{Sel}_E(\mathbb{Q}_\infty)_p = S_A(\mathbb{Q}_\infty)$. The nonprimitive Selmer groups $\operatorname{Sel}_E^{\Sigma_0}(\mathbb{Q}_\infty)_p$ and $S_A^{\Sigma_0}(\mathbb{Q}_\infty)$ also coincide, and therefore we can study the algebraic Iwasawa invariants associated to $E$ by using the results of this section. It will be useful to point out that, for an odd prime $p$, the group $S_{A[p]}(\mathbb{Q}_\infty)$ is determined just by the isomorphism class of $A[p]$ as a $G_{\mathbb{Q}}$-module. This follows from the remark at the end of the previous paragraph.

By Kato's theorem, $S_A(\mathbb{Q}_\infty)$ is $\Lambda$-cotorsion. Also, $H^0(\mathbb{Q}_\infty, A^*)$ is finite since $A^* \cong E[p^\infty]$ by the Weil pairing and $E(\mathbb{Q}_\infty)_{\text{tors}}$ is known to be finite. (If we assume that $E[p]$ is an irreducible $G_{\mathbb{Q}}$-module, then it's easy to see that $H^0(\mathbb{Q}_\infty, E[p^\infty]) = 0$.) Corollary (2.3) then implies the important relationships (6) and (7). By proposition (2.4), $\mathscr{P}_\ell(T)$ is a generator of the characteristic ideal of $\mathcal{H}_\ell(\mathbb{Q}_\infty)\widehat{\phantom{)}}$. Relationships (8) and (9) follow from this. Thus the equivalences given in theorem (1.5) have been established.

Suppose now that $E_1$ and $E_2$ are elliptic curves satisfying the hypotheses in theorem (1.4). Let $\Sigma$ be a finite set of primes containing $p$, $\infty$, and all primes where either $E_1$ or $E_2$ has bad reduction. Let $\Sigma_0 = \Sigma - \{p, \infty\}$. For $i = 1, 2$, we define $A_i = E_i[p^\infty]$. Then, by the above remarks and by proposition (2.8), we have

$$\operatorname{Sel}_{E_i}^{\Sigma_0}(\mathbb{Q}_\infty)[p] = S_{A_i}^{\Sigma_0}(\mathbb{Q}_\infty)[p] \cong S_{A_i[p]}^{\Sigma_0}(\mathbb{Q}_\infty).$$

Furthermore, the order of this group is independent of $i$ since $A_1[p] \cong A_2[p]$ as $G_{\mathbb{Q}}$-modules.

As a consequence, we see that if $\mu^{\text{alg}}_{E_1} = 0$, then $\mu^{\text{alg}}_{E_2} = 0$ and $\lambda^{\text{alg}}_{E_2, \Sigma_0} = \lambda^{\text{alg}}_{E_1, \Sigma_0}$.

**Remark (2.10)** The arguments apply with virtually no change if the Galois modules $E_1[p^\infty]$, $E_2[p^\infty]$ are replaced by similar Galois modules associated to weight 2 eigenforms whose levels are not divisible by $p^2$ and which are ordinary at $p$ (i.e., $a_p$ is a unit). Let $f_1, f_2$ be two such eigenforms. Fixing an embedding $\overline{\mathbb{Q}} \to \overline{\mathbb{Q}}_p$, the Galois module $A_i$ corresponding to $f_i$ is an $\mathcal{O}_i$-module, cofree of corank 2, where $\mathcal{O}_i$ is the closure of the integers in the field generated over $\mathbb{Q}$ by the coefficients of $f_i$. (Note that $A_i$ is defined as $V_p(f_i)/T_p(f_i)$, where $V_p(f_i)$ is the 2-dimensional $F_p$-representation associated to $f_i$. Here $F_p$ is the completion of $F$ at the prime induced by the above embedding and $T_p(f_i)$ is a $G_\mathbb{Q}$-invariant $\mathcal{O}_i$-lattice in $V_p(f_i)$.) Let $\pi_i$ denote a uniformizing parameter in $\mathcal{O}_i$. Then it is enough to assume that $A_1[\pi_1]$ and $A_2[\pi_2]$ are irreducible and become isomorphic as $G_\mathbb{Q}$-modules after extending scalars to a finite field $k$ containing both $\mathcal{O}_1/(\pi_1)$ and $\mathcal{O}_2/(\pi_2)$. In particular, we can twist the Tate modules for $E_1$ and $E_2$ satisfying the hypotheses we made above by any Dirichlet character $\chi$ whose conductor is not divisible by $p^2$ so that $\widetilde{E}_i[p^\infty] \otimes \chi$ is either unramified or tamely ramified at $p$. In some cases, remarks 2.7 and 2.9 are needed to give the desired conclusions.

The above discussion shows that if $p$ is an odd prime, $E_1[p] \cong E_2[p]$ as $G_\mathbb{Q}$-modules, and $\mu^{\text{alg}}_{E_1} = 0$, then one can compute $\lambda^{\text{alg}}_{E_2}$ if one knows $\lambda^{\text{alg}}_{E_1}$. To illustrate this, we consider the example mentioned in the introduction. The elliptic curves $E_1$ and $E_2$ have conductors 52 and 364, respectively. We take $p = 5$. One can show that $\text{Sel}_{E_1}(\mathbb{Q}_\infty)_5 = 0$. (Several examples discussed in chapter 5 of [Gre99] are completely analogous. The crucial ingredients are that $\text{Sel}_{E_1}(\mathbb{Q})_5 = 0$, that 5 is not an anomalous prime for $E_1$ (i.e., that $a_5(E_1) \not\equiv 1 \pmod 5$), and that the Tamagawa factors for $E_1$ corresponding to 2 and 13 are not divisible by 5.) It follows that $\mu^{\text{alg}}_{E_1} = \lambda^{\text{alg}}_{E_1} = 0$. One can also show that $\mu^{\text{anal}}_{E_1} = \lambda^{\text{anal}}_{E_1} = 0$. (The crucial ingredients for this are that $\mathscr{L}(E/\mathbb{Q}, T)$ is in $\Lambda$, which is known since $E_1[5]$ is irreducible, and that the interpolation property shows that $\mathscr{L}(E/\mathbb{Q}, 0) \in \mathbb{Z}_p^\times$. Hence $\mathscr{L}(E/\mathbb{Q}, T) \in \Lambda^\times$.) Let $\Sigma = \{5, \infty, 2, 7, 13\}$, $\Sigma_0 = \{2, 7, 13\}$. We need just determine $\sigma^{(\ell)}_{E_i}$ for $\ell = 2, 7, 13$, $i = 1, 2$. For $\ell = 2$, both $E_1$ and $E_2$ have additive reduction and so the corresponding Euler factor in $L(E_i/\mathbb{Q}, s)$ is 1. Thus $\sigma^{(2)}_{E_1} = \sigma^{(2)}_{E_2} = 0$. For $\ell = 13$, both $E_1$ and $E_2$ have nonsplit, multiplicative reduction. The corresponding Euler factors are both $1 + 13^{-s}$, whose value at $s = 1$ is a 5-adic unit. Thus, $\mathscr{P}_{13}(0) \in \mathbb{Z}_p^\times$ and so $\mathscr{P}_{13}(T) \in \Lambda^\times$. We have $\sigma^{(13)}_{E_1} = \sigma^{(13)}_{E_2} = 0$. However, for $\ell = 7$, $E_1$ has good reduction and $E_2$ has multiplicative reduction. The corresponding Euler factors are $1 + 2X + 7X^2$ and $1 - X$, where $X = 7^{-s}$. At $s = 1$, the value of the Euler factor for $E_2$ is in $\mathbb{Z}_p^\times$ and so $\sigma^{(7)}_{E_2} = 0$. But $1 + 2X + 7X^2 \equiv (1 - X)(1 - 2X) \pmod 5$ and $X = \widetilde{7}^{-1}$ has multiplicity 1 as a root. Also $5^2 \| (7^4 - 1)$ and so $s_7 = 5$, in the notation of proposition 2.4. That is, 7 splits completely in $\mathbb{Q}_1/\mathbb{Q}$, and the primes of $\mathbb{Q}_1$ above 7 remain inert in $\mathbb{Q}_\infty/\mathbb{Q}_1$. It follows that $\sigma^{(7)}_{E_1} = 5$. Therefore, we find that $\lambda_{E_2}^{\text{alg}} = \lambda_{E_2,\Sigma_0}^{\text{alg}} = \lambda_{E_1,\Sigma_0}^{\text{alg}} = \lambda_{E_1}^{\text{alg}} + 5 = 5$, as stated in the introduction.

Now assume that $E$ is an elliptic curve/$\mathbb{Q}$ such that $E[p]$ is reducible as a $G_{\mathbb{Q}}$-module. Assume also that $E$ has good ordinary reduction at $p$ and that $p$ is odd. Then there is an exact sequence of $G_{\mathbb{Q}}$-modules

$$0 \to \Phi \to E[p] \to \Psi \to 0$$

where both $\Phi$ and $\Psi$ are cyclic of order $p$. $G_{\mathbb{Q}}$ acts on $\Phi$ by a character $\varphi : G_{\mathbb{Q}} \to (\mathbb{Z}/p\mathbb{Z})^{\times}$, and on $\Psi$ by a character $\psi$. We can view both $\varphi$ and $\psi$ as having values in $\mathbb{Z}_p^{\times}$. Then $\varphi\psi = \omega$, the Teichmüller character of $G_{\mathbb{Q}}$. One of the characters $\varphi$ or $\psi$ is even, one odd. One of these characters is ramified at $p$, the other unramified. We will assume that the ramified character is even, and so the unramified character is odd. A result of Schneider then implies that the $\mu$-invariant for $\text{Sel}_E(\mathbb{Q}_\infty)_p^{\widehat{\phantom{p}}}$ is unchanged by a $p$-isogeny. The $\lambda$-invariant is always unchanged by an isogeny. Thus, we may assume for our purpose that $\varphi$ is ramified and even, $\psi$ is unramified and odd. It follows immediately that $H^0(\mathbb{Q}, E[p]) = 0$. Therefore by proposition 2.8 and the fact that $\text{Sel}_E^{\Sigma_0}(\mathbb{Q}_\infty)_p$ and $S_{E[p^\infty]}^{\Sigma_0}(\mathbb{Q}_\infty)$ coincide, we have

$$\text{Sel}_E^{\Sigma_0}(\mathbb{Q}_\infty)[p] \cong S_{E[p]}^{\Sigma_0}(\mathbb{Q}_\infty)$$

where $\Sigma$ is as before and $\Sigma_0 = \Sigma - \{p, \infty\}$. Now $H^0(\mathbb{Q}_\infty, \Psi) = 0$. Later we will show that $H^2(\mathbb{Q}_\Sigma/\mathbb{Q}_\infty, \Phi) = 0$. Therefore, we have an exact sequence

$$0 \to H^1(\mathbb{Q}_\Sigma/\mathbb{Q}_\infty, \Phi) \xrightarrow{\epsilon} H^1(\mathbb{Q}_\Sigma/\mathbb{Q}_\infty, E[p]) \xrightarrow{\delta} H^1(\mathbb{Q}_\Sigma/\mathbb{Q}_\infty, \Psi) \to 0.$$

Now $S_{E[p]}^{\Sigma_0}(\mathbb{Q}_\infty) = \ker\left(H^1(\mathbb{Q}_\Sigma/\mathbb{Q}_\infty, E[p]) \to H^1(I_p, \Psi)\right)$. Hence $\text{im}(\epsilon) \subset S_{E[p]}^{\Sigma_0}(\mathbb{Q}_\infty) = \delta^{-1}(\mathscr{U})$, where $\mathscr{U} = \ker\left(H^1(\mathbb{Q}_\Sigma/\mathbb{Q}_\infty, \Psi)\right) \to H^1(I_p, \Psi)$. Therefore, we have

$$\dim_{\mathbb{Z}/p\mathbb{Z}}(S_{E[p]}^{\Sigma_0}(\mathbb{Q}_\infty)) = \dim_{\mathbb{Z}/p\mathbb{Z}}(H^1(\mathbb{Q}_\Sigma/\mathbb{Q}_\infty, \Phi)) + \dim_{\mathbb{Z}/p\mathbb{Z}}(\mathscr{U}).$$

Next we will relate the two dimensions on the right to classical Iwasawa theory.

For this purpose, we consider Selmer groups for one-dimensional representations. Assume that $\dim(V_p) = 1$. Then $\text{Gal}(\mathbb{Q}_\Sigma/\mathbb{Q})$ acts on $V_p$ by a continuous homomorphism $\theta : \text{Gal}(\mathbb{Q}_\Sigma/\mathbb{Q}) \to \mathcal{O}^{\times}$. Clearly, $\theta$ factors through $G = \text{Gal}(K_\infty/\mathbb{Q})$, where $K_\infty$ is some finite extension of $\mathbb{Q}_\infty$, and $G$ is abelian. Consider the restriction map

$$H^1(\mathbb{Q}_\Sigma/\mathbb{Q}_\infty, A) \to H^1(\mathbb{Q}_\Sigma/K_\infty, A)^{\Delta}$$

where $\Delta = \text{Gal}(K_\infty/\mathbb{Q}_\infty)$. The kernel and cokernel are finite; they are trivial if $p \nmid |\Delta|$. We will assume that $p \nmid |\Delta|$, which will be sufficient for our purpose. Since $\text{Gal}(\mathbb{Q}_\Sigma/K_\infty)$ acts trivially on $A$, we have

$$H^1(\mathbb{Q}_\Sigma/K_\infty, A) = \text{Hom}(X_\infty^{\Sigma}, A)$$

where $X_\infty^\Sigma = \mathrm{Gal}\,(M_\infty^\Sigma/K_\infty)$, and $M_\infty^\Sigma$ denotes the maximal abelian, pro-$p$ extension of $K_\infty$ unramified outside $\Sigma$. We can identify $\Gamma$ with a subgroup of $G$ so that $G = \Delta \times \Gamma$. This decomposition is canonical since we are assuming that $p \nmid |\Delta|$. We have $\mathbb{Z}_p[[G]] = \mathbb{Z}_p[[\Gamma]][\Delta]$. Let $M_\infty$ denote the maximal abelian pro-$p$ extension of $K_\infty$ unramified outside $\{p, \infty\}$. Let $L_\infty$ denote the maximal abelian pro-$p$ extension of $K_\infty$ which is unramified everywhere. Depending on the parity of $\theta$, we can describe the Selmer group $S_A(\mathbb{Q}_\infty)$ in terms of either $X_\infty = \mathrm{Gal}\,(M_\infty/K_\infty)$ or $Y_\infty = \mathrm{Gal}\,(L_\infty/K_\infty)$.

Assume first that $\theta$ is even. Let $\xi = \theta|_\Delta$, which is an even character of $\Delta$. We have $d^+ = 1$ and so we must take $W_p = V_p$. Therefore, $\mathcal{H}_p(\mathbb{Q}_\infty, A) = 0$ and it follows that

$$S_A(\mathbb{Q}_\infty) \cong \mathrm{Hom}_\Delta(X_\infty, A) = \mathrm{Hom}_\mathcal{O}((X_\infty \otimes_{\mathbb{Z}_p} \mathcal{O})^\xi, A).$$

Now $K_\infty$ is the cyclotomic $\mathbb{Z}_p$-extension of $K = K_\infty^\Gamma$ and $\Delta = \mathrm{Gal}\,(K/\mathbb{Q})$ is abelian. The Ferrero-Washington theorem implies that the torsion $\Lambda$-module $(X_\infty \otimes_{\mathbb{Z}_p} \mathcal{O})^\xi$ has $\mu$-invariant equal to zero. We denote its $\mathcal{O}$-rank by $\lambda_\xi$. Then, $\lambda_\xi = \mathrm{corank}_\mathcal{O}(S_A(\mathbb{Q}_\infty))$. Note that the $\mathcal{O}$-corank of $S_A(\mathbb{Q}_\infty)$ depends only on $\xi = \theta|_\Delta$, although its structure as a $\Lambda$-module does depend on $\theta$ itself. Let $\Sigma_0 = \Sigma - \{p, \infty\}$. Let $\lambda_{\xi,\Sigma_0} = \mathrm{corank}_\mathcal{O}(S_A^{\Sigma_0}(\mathbb{Q}_\infty))$. Assuming that $\xi$ is nontrivial (so that $H^0(\mathbb{Q}, A[\pi]) = 0$), proposition (2.6) implies that $\lambda_{\xi,\Sigma_0} = \dim_{\mathcal{O}/\pi\mathcal{O}}(S_{A[\pi]}^{\Sigma_0}(\mathbb{Q}_\infty))$. Also, noting that the action of $G_\mathbb{Q}$ on $V_p^* = \mathrm{Hom}(V_p, \mathcal{O}_p(1))$ is odd, it follows that $H^0(\mathbb{Q}_\infty, A^*) = 0$. Then corollary (2.3) allows us to compute $\lambda_{\xi,\Sigma_0}$ in terms of $\lambda_\xi$.

Assume now that $\theta$ is odd. Then $\xi$ is an odd character of $\Delta$. We have $d^+ = 0$ and so we must take $W_p = 0$. That is, the local condition at $\eta_p$ occurring in the definition of $S_A(\mathbb{Q}_\infty)$ is that a cocycle class be unramified. Thus $S_A(\mathbb{Q}_\infty) = H^1_{\mathrm{unr}}(\mathbb{Q}_\Sigma/\mathbb{Q}_\infty, A)$, the group of everywhere unramified cocycle classes. It follows that

$$S_A(\mathbb{Q}_\infty) = \mathrm{Hom}_\Delta(Y_\infty, A) = \mathrm{Hom}((Y_\infty \otimes_{\mathbb{Z}_p} \mathcal{O})^\xi, A).$$

Again, the Ferrero-Washington theorem implies that the $\Lambda$-torsion module $(Y_\infty \otimes_{\mathbb{Z}_p} \mathcal{O})^\xi$ has $\mu$-invariant zero. We let $\lambda_\xi$ denote its $\mathcal{O}$-rank, and so $S_A(\mathbb{Q}_\infty)$ has $\mathcal{O}$-corank equal to $\lambda_\xi$. As we remarked above, this clearly depends only on $\xi = \theta|_\Delta$. Assume that $\xi \ne \omega$. Then $H^0(\mathbb{Q}_\infty, A^*) = 0$ and so we can apply corollary (2.3) to determine $\lambda_{\xi,\Sigma_0} = \mathrm{corank}_\mathcal{O}(S_A^{\Sigma_0}(\mathbb{Q}_\infty))$ in terms of $\lambda_\xi$. Since $\xi$ is odd, we have $H^0(\mathbb{Q}, A[\pi]) = 0$ and therefore proposition (2.8) implies that $\lambda_{\xi,\Sigma_0} = \dim_{\mathcal{O}/\pi\mathcal{O}}(S_{A[\pi]}^{\Sigma_0}(\mathcal{O}_\infty))$.

We can apply these observations to $\theta = \varphi$ and $\theta = \psi$, regarding $\Phi$ or $\Psi$ as $A[p]$ where $A \cong \mathbb{Q}_p/\mathbb{Z}_p$ is a group and $G_\mathbb{Q}$ acts by either $\varphi$ or $\psi$. Now $\varphi$ is even, but it is ramified and hence nontrivial. Also, $\psi$ is odd but since $\psi = \omega\varphi^{-1}$, we have $\psi \ne \omega$. Hence we obtain the result that $\mu_E^{\mathrm{alg}} = \mu_{E,\Sigma_0}^{\mathrm{alg}} = 0$ and that

$$\lambda_{E,\Sigma_0}^{\mathrm{alg}} = \lambda_{\varphi,\Sigma_0} + \lambda_{\psi,\Sigma_0} \tag{16}$$

which shows that $\lambda^{\text{alg}}_{E,\Sigma_0}$ can be easily calculated in terms of the classical Iwasawa invariant $\lambda_\varphi = \lambda_\psi$. We must just explain why $H^2(\mathbb{Q}_\Sigma/\mathbb{Q}_\infty, \Phi) = 0$. This is in fact equivalent to the vanishing of the $\mu$-invariant for the $\Lambda$-module $X^\varphi_\infty$. We take $\theta = \varphi$. Suppose that $A \cong \mathbb{Q}_p/\mathbb{Z}_p$ and that $G_\mathbb{Q}$ acts on $A$ by $\varphi$. Since $X^\varphi_\infty$ is $\Lambda$-torsion (because $\varphi$ is even), it follows that $S_A(\mathbb{Q}_\infty)$ and $S^{\Sigma_0}_A(\mathbb{Q}_\infty) = H^1(\mathbb{Q}_\Sigma/\mathbb{Q}_\infty, A)$ are $\Lambda$-cotorsion. Their $\mu$-invariants are zero, by the Ferrero-Washington theorem. By proposition (2.5), $S^{\Sigma_0}_A(\mathbb{Q}_\infty)\widehat{\phantom{x}}$ has no nonzero, finite $\Lambda$-submodules, which then implies that it is a divisible group. That is, $H^1(\mathbb{Q}_\Sigma/\mathbb{Q}_\infty, A)$ is divisible. Also, $H^2(\mathbb{Q}_\Sigma/\mathbb{Q}_\infty, A)$ must be $\Lambda$-cotorsion (using the first result about coranks recalled at the beginning of this section). By proposition 4 of [Gre99], it follows that $H^2(\mathbb{Q}_\Sigma/\mathbb{Q}_\infty, A) = 0$. Then the exact sequence $0 \to \Phi \to A \xrightarrow{p} A \to 0$ induces an isomorphism

$$H^1(\mathbb{Q}_\Sigma/\mathbb{Q}_\infty, A)/pH^1(\mathbb{Q}_\Sigma/\mathbb{Q}_\infty, A) \xrightarrow{\sim} H^2(\mathbb{Q}_\Sigma/\mathbb{Q}_\infty, \Phi).$$

The vanishing of $H^2(\mathbb{Q}_\Sigma/\mathbb{Q}_\infty, \Phi)$ follows from this.

In the above discussion we may assume that $K_\infty$ contains $\mu_p$, and hence $\mu_{p^\infty}$. (It is not necessary for $\theta$ to be a faithful character of $\text{Gal}(K_\infty/K)$.) Now complex conjugation (in $\text{Gal}(K_\infty/\mathbb{Q})$) acts on both $X_\infty$ and $Y_\infty$, and one can then define $X^\pm_\infty$, $Y^\pm_\infty$ as the ($\pm 1$)-eigenspaces for complex conjugation. A theorem of Iwasawa then implies that $Y^-_\infty$ is pseudo-isomorphic to $\text{Hom}(X^+_\infty, \mathbb{Z}_p(1))$ as a module for $\mathbb{Z}_p[[G]] = \Lambda[\Delta]$, where $\mathbb{Z}_p(1)$ just denotes $\mathbb{Z}_p$ together with an action of $G$ by the $p$-cyclotomic character $\chi$. (Stated this way, the result also depends on the vanishing of $\mu$ for $X^+_\infty$ and $Y^-_\infty$.) It follows that, if $\xi$ is an even character of $\Delta$, then $(X_\infty \otimes_{\mathbb{Z}_p} \mathcal{O})^\xi$ has the same $\mathcal{O}$-rank as $(Y_\infty \otimes_{\mathbb{Z}_p} \mathcal{O})^{\omega \xi^{-1}}$, and therefore $\lambda_\xi = \lambda_{\omega \xi^{-1}}$. In particular, $\lambda_\varphi = \lambda_\psi$. The easiest way to compute $\lambda_\varphi$ is to use the theorem of Mazur and Wiles which implies that $\lambda_\varphi$ is just the $\lambda$-invariant of the Kubota-Leopoldt $p$-adic $L$-function $L_p(\varphi, s)$. Extensive calculations have been carried out by T. Fukuda when $\psi$ is imaginary quadratic and $p = 3, 5$, or $7$.

For a character $\theta \colon \text{Gal}(\mathbb{Q}_\Sigma/\mathbb{Q}) \to \mathbb{Z}_p^\times$, let $A_\theta$ denote the group $\mathbb{Q}_p/\mathbb{Z}_p$ on which $\text{Gal}(\mathbb{Q}_\Sigma/\mathbb{Q})$ acts by $\theta$. Recall that for any prime $l \neq p$, $s_l$ denotes the number of primes $\eta$ of $\mathbb{Q}_\infty$ lying above $l$, or equivalently the $p$-adic valuation of $(l^{p-1} - 1)/p$. Define $t_l(E)$ as the integer

$$\text{corank}_{\mathbb{Z}_p}(H^1((\mathbb{Q}_\infty)_\eta, A_\varphi)) + \text{corank}_{\mathbb{Z}_p}(H^1((\mathbb{Q}_\infty)_\eta, A_\psi)) - \text{corank}_{\mathbb{Z}_p}(H^1((\mathbb{Q}_\infty)_\eta, E[p^\infty]))$$

where $\eta$ is any prime of $\mathbb{Q}_\infty$ lying over $l$. As a consequence of corollary (2.3) applied to $E[p^\infty]$, $A_\varphi$ and $A_\psi$, together with (16), we see that

$$\lambda_E = 2\lambda_\psi + \sum_{l \in \Sigma_0} s_l t_l(E).$$

One way to calculate the above $\mathbb{Z}_p$-coranks is to use proposition 2 of [Gre89]. If $A$ is a $G_{\mathbb{Q}_l}$-module (where $l \neq p$) which is isomorphic to $(\mathbb{Q}_p/\mathbb{Z}_p)^d$ as a group, let $T_A$ denote its Tate module. Define $A^* = \text{Hom}(T_A, \mu_{p^\infty})$, which is also a $G_{\mathbb{Q}_l}$-module isomorphic to $(\mathbb{Q}_p/\mathbb{Z}_p)^d$. Then the above cited proposition implies that $H^1((\mathbb{Q}_\infty)_\eta, A)$ has the same $\mathbb{Z}_p$-corank as $H^0((\mathbb{Q}_\infty)_\eta, A^*)$, which is denoted more briefly by $A^*((\mathbb{Q}_\infty)_\eta)$ in [Gre89]. The Weil pairing implies that $E[p^\infty]^* \cong E[p^\infty]$ as $G_{\mathbb{Q}_l}$-modules. Also, since $\varphi\psi = \omega$ and $\chi|_\Delta = \omega$, where $\chi$ is the $p$-cyclotomic character as before, it is clear that $A_\varphi^* \cong A_\psi$ and $A_\psi^* \cong A_\varphi$ as $G_{(\mathbb{Q}_\infty)_\eta}$-modules. Thus

$$t_l(E) = \text{corank}_{\mathbb{Z}_p}(A_\varphi((\mathbb{Q}_\infty)_\eta)) + \text{corank}_{\mathbb{Z}_p}(A_\psi((\mathbb{Q}_\infty)_\eta)) - \text{corank}_{\mathbb{Z}_p}(E((\mathbb{Q}_\infty)_\eta)_{p\text{-tors}}).$$

We can now discuss the example mentioned in the introduction where $J$ is one of the curves of conductor 11, $E$ is the quadratic twist $J_{-c}$, and $p = 5$. The character $\psi$ is the quadratic character corresponding to $F = \mathbb{Q}(\sqrt{-c})$. Let $\Sigma$ consist of $\infty$, 5, 11, and all the ramified primes for $F/\mathbb{Q}$. Let $\Sigma_0 = \Sigma - \{\infty, 5\}$. Assume first that $l$ is ramified in $F/\mathbb{Q}$. Then both $\varphi = \omega\psi^{-1}$ and $\psi$ are nontrivial characters of $G_{\mathbb{Q}_l}$ with orders not divisible by $p$. Their restrictions to $G_{(\mathbb{Q}_\infty)_\eta}$ are also nontrivial and so $A_\varphi((\mathbb{Q}_\infty)_\eta)$, $A_\psi((\mathbb{Q}_\infty)_\eta)$, and $E((\mathbb{Q}_\infty)_\eta)_{p\text{-tors}}$ are all trivial. Hence $t_l(E) = 0$ for those $l$'s. This also applies to $l = 11$ if $\psi|_{G_{\mathbb{Q}_{11}}}$ is nontrivial, i.e., if 11 is ramified or inert in $F/\mathbb{Q}$. (Note that then $\varphi|_{G_{\mathbb{Q}_{11}}}$ is also nontrivial since $\omega|_{G_{\mathbb{Q}_{11}}}$ is trivial.) But if 11 splits in $F/\mathbb{Q}$, then both $\varphi|_{G_{\mathbb{Q}_{11}}}$ and $\psi_{G_{\mathbb{Q}_{11}}}$ are trivial characters. Hence $A_\varphi((\mathbb{Q}_\infty)_\eta)$ and $A_\psi((\mathbb{Q}_\infty)_\eta)$ both have $\mathbb{Z}_p$-corank 1. Now $E \cong J$ over $\mathbb{Q}_{11}$ and has split, multiplicative reduction. Since $\mu_p \subset \mathbb{Q}_{11}$, we have $(\mathbb{Q}_\infty)_\eta = \mathbb{Q}_{11}(\mu_{p^\infty})$ where $\eta$ is the unique prime of $\mathbb{Q}_\infty$ lying above 11. (Note that $s_{11} = 1$.) We have $\text{corank}_{\mathbb{Z}_p}(E((\mathbb{Q}_\infty)_\eta)_{p\text{-tors}}) = 1$. Thus $t_{11}(E) = 1$ if 11 splits in $F/\mathbb{Q}$. These remarks imply that $\lambda_E^{\text{alg}} = 2\lambda_\psi + \epsilon_\psi$ as stated in the introduction.

## 3 Congruences for $p$-adic L-functions.

In this section we prove the congruences for $p$-adic L-functions described in the introduction. The technical tools for this were developed in [Vat97], where certain canonical periods, well suited to the study of congruences, were attached to a cuspform $f$. We want to specialize this to the case where $f$ has rational Fourier coefficients and corresponds to a modular elliptic curve. Our main task in this section is therefore to compare the canonical periods of $f$ to the Néron periods of the associated elliptic curve.

At the end of this section we give a brief discussion of our results as they apply to modular forms whose Fourier coefficients are not necessarily rational numbers (so the associated abelian varieties need not be elliptic curves). While much of the theory goes through without change, there are some serious complications. In the first place, we want to apply the congruence techniques of [Vat97], and to do this we need to assume that the maximal ideals arising in the Hecke ring satisfy a certain technical condition (we assume that they are $p$-distinguished; see below for the precise definition). Secondly, we would like to compare the periods of the modular forms with those of the associated abelian varieties, and this can only be achieved under a certain semistability assumption. Neither of these assumptions need be true in general for forms on $\Gamma_1(N)$ if $p|N$. Both assumptions will be true, however, if $(N, p) = 1$, or if $f$ corresponds to an elliptic curve. For simplicity of exposition, we have chosen to concentrate on the elliptic curve case in this section. We have included statements of the general results, and the interested reader should have no difficulty filling in the details.

## Canonical periods

We begin by recalling the construction of [Vat97]. Let $E$ be a modular elliptic curve of conductor $N$, and let $p$ be a fixed odd prime. We assume that $E$ has either good ordinary or multiplicative reduction at $p$, corresponding to the two cases $(N, p) = 1$ and $(p, N/p) = 1$ respectively. We will write $f = f(z)$ for the newform of level $N$ corresponding to $E$. Then $f$ has rational Fourier coefficients. Let $\Gamma$ denote the group $\Gamma_1(N)$; we may assume that $N > 4$, so that $\Gamma$ is torsion-free. Let $S_2(\Gamma, \mathbb{Z}_p)$ denote the space of cuspforms of weight 2 on $\Gamma$ with coefficients in $\mathbb{Z}_p$, and write $H^1_{\mathrm{par}}(\Gamma, \mathbb{Z}_p)$ for the parabolic cohomology group of Eichler-Shimura. There is a natural action of complex conjugation on this cohomology group, and we denote the $(\pm 1)$-eigenspaces by a superscript. Let $\mathfrak{m}$ denote the maximal ideal above $p$ cut out by $f$ in the Hecke ring $\mathbf{T}$ generated by the usual operators $T_q, U_\ell, <q>$ over $\mathbb{Z}$. Let $\mathbf{T}_\mathfrak{m}$ denote the completion of $\mathbf{T}$ at $\mathfrak{m}$. Then $\mathbf{T}_\mathfrak{m}$ is a finite flat $\mathbb{Z}_p$-algebra. The Hecke operators act on this cohomology group, and on the space of cuspforms with integral coefficients.

There is a semisimple representation $\rho_\mathfrak{m} : G_\mathbb{Q} \to GL_2(\mathbf{T}/\mathfrak{m}) = GL_2(\mathbb{F}_p)$ characterized by $\mathrm{Tr}\,(\mathrm{Frob}(q)) = T_q$ for $(q, Np) = 1$; this representation $\rho_\mathfrak{m}$ is the semisimplification of the representation of $G_\mathbb{Q}$ on $E[p]$. If $E[p]$ is irreducible, then $\rho_\mathfrak{m}$ is equivalent to the $G_\mathbb{Q}$-representation on $E[p]$. If $E[p]$ is reducible, then the Weil pairing implies that there exists a character $\psi : \mathrm{Gal}\,(\overline{\mathbb{Q}}/\mathbb{Q}) \to \mathbb{F}_p^\times$ such that the Jordan-Holder factors of $E[p]$ are given by $\psi$ and $\omega\psi^{-1}$. Here $\omega : \mathrm{Gal}\,(\overline{\mathbb{Q}}/\mathbb{Q}) \to \mathbb{F}_p^\times$ denotes the Teichmüller character. If we let $\varphi = \omega\psi^{-1}$, then the semi-simple representation $\rho_\mathfrak{m}$ is given explicitly as

$$\rho_\mathfrak{m} = \varphi \oplus \psi. \tag{17}$$

We may assume that $\psi$ is unramified at $p$. Then we choose a sign $\alpha = \pm$ as follows:

1. If $E[p]$ is irreducible, then $\alpha$ is arbitrary.

2. If $E[p]$ is reducible, then $\alpha$ is determined by $\pm 1 = -\psi(-1)$.

A choice of sign satisfying one of these conditions is said to be *admissible* (for $E$ and $p$). Now let $f = \sum a_n q^n$ denote the newform associated to $E$. Let $M$ be any integer divisible by $N$, and let $g = \sum b_n q^n$ denote any modular form of level $M$ that is an eigenform for the *full* Hecke algebra $\mathbf{T}_1(M)$ for $\Gamma_1(M)$, and which satisfies $a_n = b_n$ whenever $(n, M) = 1$. If $\alpha$ is admissible, it follows from work of Mazur, Ribet, Wiles, and others (see [Vat97], Theorems 1.3 and 2.7) that there is an isomorphism of $\mathbf{T}_1(M)_\mathfrak{m}$–modules

$$\theta^\alpha : S_2(\Gamma_1(M), \mathbb{Z}_p)_\mathfrak{m} \cong H^1_{\mathrm{par}}(\Gamma_1(M), \mathbb{Z}_p)^\alpha_\mathfrak{m}. \tag{18}$$

Here $\mathfrak{m}$ denotes the maximal ideal of $\mathbf{T}_1(M)$ induced by the homomorphism $\pi_g : \mathbf{T}_1(M) \to \mathcal{O}$, and $\mathbf{T}_1(M)_\mathfrak{m}$ denotes the completion. Thus, if $\alpha$ is admissible, we may define canonical cocycles $\delta_g^\alpha = \theta^\alpha(g) \in H^1_{\mathrm{par}}(\Gamma_1(M), \mathbb{Z}_p)$. Now consider the differential form $g(z)dz$ on the upper-half-plane. Then we define a cocycle $\omega_g \in H^1_{\mathrm{par}}(\Gamma_1(M), \mathbb{C})$ by the usual Eichler-Shimura construction: if $\gamma \in \Gamma_1(M)$, then

$$\omega_g(\gamma) = \int_{z_0}^{\gamma z_0} g(z)dz.$$

We may decompose $\omega_g$ into plus and minus parts $\omega_g^\pm$ according to the action of complex conjugation. Note that each of the cocycles $\delta_g^\alpha$ and $\omega_g^\pm$ is an eigenvector for the action of the *full* Hecke algebra, since the form $g$ is so. In each case, the eigenvalue for $T_q$ is $a_q$, when $q$ is prime to $M$, and the eigenvalue for $U_q$ coincides with the eigenvalue of $U_q$ on $g$. Thus we find that there exist complex periods $\Omega_g^\alpha$ such that

$$\Omega_g^\alpha \delta_g^\alpha = \omega_g^\alpha. \tag{19}$$

The numbers $\Omega_g^\alpha$ are the so-called "canonical periods" for $g$ (when they exist).

Now we return to the elliptic curve $E$. For each admissible sign $\alpha$, we let $c_\alpha$ denote a generator of the $\mathbb{Z}$-module $H_1(E, \mathbb{Z})^\alpha \cong \mathbb{Z}$. Let $\omega_E$ denote a Néron differential on $E$; then we define a Néron period for $E$ by

$$\Omega_E^\alpha = \int_{c_\alpha} \omega_E. \tag{20}$$

We want to relate this Néron period on $E$ to the canonical period of $f$ by studying the geometry of a modular parametrization $X_1(N) \to E$. To do this we need to specify a model of $X_1(N)$ over $\mathbb{Q}$. For the present purpose, it is convenient to use the conventions of [DI95], Variant 9.3.6 (see also [Ste89], p. 80). Thus, $X_1(N)_\mathbb{Q}$ denotes the scheme classifying generalized elliptic curves $A$ together with an embedding of group schemes $\mu_N \hookrightarrow A$. With this convention, the cusp infinity is rational over $\mathbb{Q}$. Given a $\mathbb{Q}$-isogeny class $\mathcal{A}$ of elliptic curves, we call $E \in \mathcal{A}$ *optimal* if there exists a modular parametrization

$$\pi : X_1(N)_{\mathbb{Q}} \to E$$

such that the following equivalent conditions are met:

1. The induced map on homology $H_1(X_1(N), \mathbb{Z}) \to H_1(E, \mathbb{Z})$ is surjective.

2. The induced map $J_1(N) \to E$ has connected kernel, so that there is an exact sequence of abelian varieties over $\mathbb{Q}$:

$$0 \to A \xrightarrow{i} J_1(N) \xrightarrow{\pi} E \to 0. \tag{21}$$

3. There exists an embedding $E \hookrightarrow J_1(N)$.

4. If $E' \in \mathcal{A}$, then any morphism $X_1(N) \to E'$ factors as $X_1(N) \xrightarrow{\pi} E \to E'$.

In this case, the map $\pi$ is said to be an optimal parametrization. Thus the optimal curve $E$ is an analogue for $X_1(N)$-parametrizations of the strong Weil curve arising from parametrizations by $X_0(N)$. Note, however, that the strong Weil curve and the optimal curve are *not* in general equal. The optimal curves were introduced and studied by Stevens [Ste89]. Then we have the following proposition:

**Proposition (3.1)** *Assume that $E$ is optimal in its isogeny class, and that $p$ is a prime of either good ordinary or multiplicative reduction for $E$. Then the numbers $\Omega_E^\alpha$ and $(-2\pi i)\Omega_f^\alpha$ are equal up to a factor which is a $p$-adic unit.*

The next several paragraphs are devoted to a proof of this assertion. The main ingredient is a study of the "Manin constant" associated to the optimal parametrization $X_1(N) \to E$, especially when $p$ is a prime of multiplicative reduction for $E$. Let $\omega_E$ denote a Néron differential on the Néron model $E_{\mathbb{Z}}$. Then $\pi^*(\omega_E) = c_1 f(q) dq/q$, for a newform $f$ on $\Gamma_1(N)$, $q = e^{2\pi i z}$, and a quantity $c_1 = c_1(\pi) \in \mathbb{Z}$ (see [Ste89], Theorem 1.6). We want to prove that $c_1$ is a $p$-adic unit, for all odd primes $p$ of either good or multiplicative reduction. The analogous result for optimal parametrizations by $X_0(N)$ is due to Mazur:

**Theorem (3.2) ([Maz78], Cor. 4.1)** *Let $E'$ be a strong Weil curve, and $p$ an odd prime of good or multiplicative reduction. Let $\pi_0 : X_0(N) \to E'$ be the strong parametrization, and $c_0$ the associated Manin constant, defined by $\pi_0^*(\omega_E) = c_0 f(q) dq/q$. Then the number $c_0 \in \mathbb{Z}$ is a $p$-adic unit.*

The proof of this theorem relies on the fact that the Jacobian $J_0(N)$ is a semi-abelian scheme over $\mathbb{Z}_p$ when $N$ is divisible by precisely the first power of $p$. Unfortunately, $J_1(N)$ is not generally semi-abelian, and Mazur's arguments cannot be transferred directly when the level $N$ is divisible by $p$. Nevertheless, we can still offer the following result.

**Proposition (3.3)** *Assume that $N$ is divisible by at most the first power of $p$. Then the integers $c_0$ and $c_1$ are equal up to a unit factor in $\mathbb{Z}_p$. The integer $c_1$ is a $p$-adic unit.*

*Proof.* In view of Mazur's theorem above, it suffices to verify the first assertion. Consider the composite

$$X_1(N)_{\mathbb{Q}} \to X_0(N)_{\mathbb{Q}} \to E'_{\mathbb{Q}} \tag{22}$$

where the first map is the natural projection induced by $\Gamma_1(N) \subset \Gamma_0(N)$, and the second is the strong parametrization. By definition of optimality, the map (22) factors as

$$X_1(N)_{\mathbb{Q}} \to E_{\mathbb{Q}} \to E'_{\mathbb{Q}}, \tag{23}$$

where the first map is optimal for $X_1$. Observe now that, under the map (22), a Néron differential on $E'$ over $\mathbb{Z}$ pulls back on $X_1(N)$ to the modular form $c_0 f(z)dz$. This follows from the definition of $c_0$, together with the fact that the natural projection of $X_1$ to $X_0$ induces the identity on $q$-expansions at infinity.

On the other hand, we can compute the pullback of a Néron differential on $E'$ by using (23) instead. The Manin constant of the first map is $c_1$, as that map is optimal for $X_1$. So the proposition will follow if we can check that a Néron differential on $E'_{\mathbb{Z}_p}$ pulls back to a Néron differential on $E_{\mathbb{Z}_p}$, for each odd prime $p$ of good or multiplicative reduction (here, as elsewhere, the subscript $\mathbb{Z}_p$ on an abelian variety over $\mathbb{Q}$ denotes the Néron model). To check this condition on differentials, it suffices to show that $E_{\mathbb{Z}_p} \to E'_{\mathbb{Z}_p}$ is étale. We will only treat the case of multiplicative reduction, as the case of good reduction is similar but easier.

We start by studying the situation over the generic fibre $\mathbb{Q}_p$. Consider the diagrams gotten from (22) and (23) by replacing $X_0$ and $X_1$ by their respective Jacobians $J_0$ and $J_1$. We get a commutative diagram over $\mathbb{Q}$:

$$\begin{array}{ccc} E' & \to & J_0 \\ \downarrow & & \downarrow \\ E & \to & J_1 \end{array}$$

where the arrows of the square are induced by Picard functoriality. The maps $E' \to J_0$ and $E \to J_1$ are both injective, and the kernel of $J_0 \to J_1$ is the so-called Shimura subgroup.

This is a finite flat group scheme (still over $\mathbb{Q}_p$) which is the Cartier dual of a constant group scheme. We see easily from this that the kernel $G$ of the dual map $E \to E'$ is constant over $\mathbb{Q}_p$. Observe that this dual map is precisely the morphism appearing in (23).

We may assume, using the Sylow theorems and the fact that multiplication by $\ell$ is finite flat and étale over $\mathbb{Z}_p$, that $G$ is nontrivial and has $p$-power order. We have shown that $G$ is constant over the generic fibre, and it will suffice to transfer the constancy from $\mathbb{Q}_p$ to $\mathbb{Z}_p$. Here we will use the fact that the Neron models of $E'$ and $E$ are both semi-abelian schemes (since the curves have multiplicative reduction). Because of the semi-abelianness, we find that the kernel $K_{\mathbb{Z}_p}$ of $E_{\mathbb{Z}_p} \to E'_{\mathbb{Z}_p}$ is quasi-finite and flat over $\mathbb{Z}_p$ (see [BRL90], pages 177-178).

Now write $G_{\mathbb{Z}_p}$ for the constant finite flat group scheme over $\mathbb{Z}_p$ whose generic fibre was called $G$ above. Then by the Néron property there is a map $G_{\mathbb{Z}_p} \to E_{\mathbb{Z}_p}$. Composing with $E_{\mathbb{Z}_p} \to E'_{\mathbb{Z}_p}$, we find that the map $G_{\mathbb{Z}_p} \to E_{\mathbb{Z}_p} \to E'_{\mathbb{Z}_p}$ is trivial along the generic fibre, hence trivial (everything in sight is flat). We find that $G$ maps into the kernel $K_{\mathbb{Z}_p}$ of $E_{\mathbb{Z}_p} \to E'_{\mathbb{Z}_p}$, and that the isomorphism $G = K$ holds over $\mathbb{Q}_p$.

Now, since $K_{\mathbb{Z}_p}$ is quasi-finite and flat, and since $\mathbb{Z}_p$ is Henselian, we have a canonical subgroup scheme $FK_{\mathbb{Z}_p} \subset K_{\mathbb{Z}_p}$, where $FK_{\mathbb{Z}_p}$ is finite flat, and $K/FK$ is quasi-finite and étale, and has trivial special fibre. Thus $G_{\mathbb{Z}_p}$ sits inside the finite part of $K_{\mathbb{Z}_p}$. Since $(K/G)_{\mathbb{Z}_p}$ is also quasi-finite and flat, and has trivial generic fibre, it follows that $(K/G)_{\mathbb{Z}_p}$ is finite flat (the quasi-finite étale part has to be trivial on both fibres, hence trivial). But $G_{\mathbb{Z}_p}$ and $K_{\mathbb{Z}_p}$ are equal along the generic fibre, and must therefore have $K = G$. This shows that $K$ is finite flat and constant, and hence that $E_{\mathbb{Z}_p} \to E'_{\mathbb{Z}_p}$ is étale, as claimed.

We can now complete the proof of (3.1). Consider the integration map

$$H_1(X_1(N), \mathbb{Z}_p)^\alpha \xrightarrow{\int f(z)dz} \mathbb{C}.$$

The image $H(f)^\alpha$ is a free $\mathbb{Z}_p$-module of rank 1, generated by the number $\Omega_f^\alpha$. To verify this last statement, note that since $a_1(f) = 1$, the form $f$ is not divisible by $p$ in the space of cuspforms with integral coefficients. Thus the cocycle $\delta_f^\alpha$ is not divisible by $p$ in the integral cohomology group, and there exists a homology class $\gamma$ such that $\delta_f^\alpha$ takes on a $p$-adic unit value when capped against $\gamma$. Since $\omega_f^\alpha = \Omega^\alpha \delta_f^\alpha$, it follows that $H(f)^\alpha$ is generated by the period $\Omega^\alpha$. If $E$ is optimal, we define a free $\mathbb{Z}_p$-module $H(E)^\alpha$ by integrating the Néron differential $\omega_E$ on cycles in $H_1(E, \mathbb{Z}_p)^\alpha$. Obviously $H(E)^\alpha$ is generated by $\Omega_E^\alpha$. Thus it suffices to show that $(2\pi i)H(f)^\alpha = H(E)^\alpha$. But this is clear, in view of the adjointness formula for the map $\pi : X_1(N) \to E$, the surjectivity of the map on homology, and the fact that the constant $c(\pi)$ is a $p$-adic unit. ∎

**Remark (3.4)** If $E$ does not admit any $p$-isogenies, so that $E[p]$ is irreducible, then it is clear that the Néron periods of any isogeneous curve differ from those of $E$ by a $p$-adic unit. Thus in this case we obtain $\Omega_E^\alpha = \Omega_f^\alpha$, for any choice of sign.

**Remark (3.5)** The proof of (3.1) also yields some information about the periods of $\Omega_E^\alpha$ when $\alpha$ is *not* admissible. Namely, it is clear from the proof that, in this case, the cocycle $\delta_f^\alpha = \omega_f^\alpha / \Omega_E^\alpha$ lies in the integral cohomology group $H^1_{\text{par}}(\Gamma, \mathbb{Z}_p)$, and that $p$ does not divide $\delta_f^\alpha$ in $H^1_{\text{par}}(\Gamma, \mathbb{Z}_p)$.

As we have already remarked in the introduction, the congruence formulae are satisfied only by the nonprimitive $p$-adic L-functions. It will therefore be useful later to have a comparison between the periods of an nonprimitive form with those of the associated newform. The necessary result is given below.

**Lemma (3.6)** *Let $f$ be a newform of level $N = \sum a_n q^n$, corresponding to the modular elliptic curve $E_{\mathbb{Q}}$. Assume that $E$ has ordinary reduction at $p$, and let $g = \sum b_n q^n$ be the eigenform, of level $M$, obtaining from $f$ by removing all Euler factors at primes $q \neq p$ such that $q|N$. Then, if $\alpha$ is an admissible sign for $E$, the canonical periods $\Omega_f^\alpha$ and $\Omega_g^\alpha$ are equal up to $p$-adic unit.*

*Proof.* Since $f$ corresponds to an elliptic curve, we have $a_n \in \mathbb{Q}$, for all $n$. Then, by construction of $g$, we have $b_n \in \mathbb{Q}$ for all $n$. Now let $B$ denote the abelian variety quotient of $J_1(M)$ associated by Shimura to $g$ in [Shi73]. Note that the existence of this abelian variety does *not* require that $g$ be a newform, merely that it be an eigenvector for all the Hecke operators. Since $g$ has rational Fourier coefficients, it is clear that $B$ is an elliptic curve isogenous to $E$. Shimura has shown furthermore that if $\omega$ denotes a Néron differential on $B$, then the pullback of $\omega$ to $J_1(M)$ under the natural quotient map $J_1(M) \to B$ is an eigenvector for the Hecke algebra, with eigenvalues equal to those of the cuspform $g$. For all these results of Shimura, we refer the reader to [Shi73], Theorem 1, page 526. We may assume that the kernel $K$ of the quotient map is connected. We contend then that the abelian variety $K$ is stable under the action of the Hecke algebra. To see this, consider the exact sequence $0 \to K \xrightarrow{j} J_1(M) \xrightarrow{\hat{\pi}} B \to 0$. Let $t$ denote the endomorphism of $J_1(M)$ induced by any Hecke operator over $\mathbb{Q}$. It suffices to show that the composite $j_t : K \to B$ given by $j_t = j \circ t \circ \hat{\pi}$ is trivial (we are viewing operators as acting on the right). Let $\omega$ denote a Néron differential on $B$, and consider the differential form $j_t^* \omega$. Then one has $(j \circ t \circ \hat{\pi})^* \omega = j^*((t \circ \hat{\pi})^* \omega) = j^*(t(g)\hat{\pi}^* \omega)$, since $\hat{\pi}^* \omega$ is an eigenvector for the Hecke algebra, with eigenvalues given by those of $g$. But now $j^*(t(g)\hat{\pi}^* \omega) = t(g) \cdot (j^* \hat{\pi}^* \omega) = t(g) \cdot (j \circ \hat{\pi})^* \omega = 0$. Thus $j_t$ induces the zero map on cotangent spaces. Since $K$ is connected, it follows that $j_t = 0$.

We will show that the canonical periods of $g$ coincide with the Néron periods of $B$. To do this, it will be enough to show that the Manin constant of the parametrization $\pi_M : X_1(M) \hookrightarrow J_1(M) \to B$ is a $p$-adic unit. First consider the case of good reduction, that is, the case that $(p, N) = (p, M) = 1$. In this case, one checks that the original argument of Mazur from (3.2) can be modified without difficulty. The only step in which any property of newforms was used was in verifying the stability of $K$ under the Hecke operators, and we have checked this above. However, the case of multiplicative reduction is more delicate. To treat this case, let $\Gamma$ denote the group $\Gamma_1(M/p) \cap \Gamma_0(p)$, and let $X$ denote the corresponding modular curve. Then $g$ is a Hecke eigenform for $\Gamma$. Let $C$ denote the elliptic curve quotient of $J = \text{Jac}(X)$, corresponding to $g$, constructed as above. Let $\pi : X \to C$ denote the corresponding parametrization, so that if $\omega$ denotes a Néron differential on $C$, then $\pi^*\omega = c \cdot g(z)dz$, for a nonzero $c \in \mathbb{Z}$. Since $X$ has semistable reduction at $p$ (see [Wil95], page 485, or [MW84], Chapter 2), Mazur's argument in [Maz78] implies that $c$ is a $p$-adic unit. Indeed, Mazur's argument shows that $\pi^*\omega$ restricts to a nonzero differential on $X_{\mathbb{F}_p}$. Now, if $p|c$, then the $q$-expansion of $\pi^*\omega$ must vanish on the component of $X_{\mathbb{F}_p}$ containing the cusp $\infty$. But Wiles has shown (See [Wil95], Lemma 2.2, especially the bottom paragraph on page 486, and recall that we are in the ordinary case.) that any nonzero differential on $X_{\mathbb{F}_p}$ is nonzero on the component containing $\infty$.[1] This implies that $c$ is a $p$-adic unit.

To sum up, we have shown that the canonical periods of the oldform $g$ coincide with the Néron periods of an elliptic curve $B$ isogenous to the optimal curve $E$ of level $N$. If $E[p]$ is irreducible as a Galois-module, then our lemma follows from Proposition (3.1), as $E$ does not admit any nontrivial $p$-isogenies. It remains therefore to treat the case that $E[p]$ is reducible. In this case, we let $E^{\min}$ denote the *minimal* curve in the isogeny class $\mathcal{A}$ of $E$ and $B$ constructed by Stevens in [Ste89], section 2. Stevens has shown that, if $A \in \mathcal{A}$, there exists an étale isogeny $\varphi : E^{\min}_{\mathbb{Z}} \to A_{\mathbb{Z}}$. If $\alpha$ is an admissible sign, it follows from the definitions that the kernel of $\varphi$ has parity $-\alpha$ for the action of complex conjugation. This implies that the periods $\Omega_A^\alpha$ and $\Omega_{E^{\min}}^\alpha$ coincide. Since both $E$ and $B$ are members of $\mathcal{A}$, the assertion of the lemma follows. ∎

## $p$-adic L-functions

We want to give the definition of the $p$-adic L-function of a modular elliptic curve $E$, together with its various twists. More generally, we will define the $p$-adic L-function of a weight-two modular form. These functions were constructed by Mazur, Tate, and Teitelbaum in [MTT86]. Thus let $K$ be an abelian number field. We assume that $K$ is unramified at all primes dividing the level $N$, and tamely ramified at $p$. The curve $E$ is assumed of course to have ordinary reduction at $p$. Put $G = \text{Gal}(K/\mathbb{Q})$, and fix a character $\chi$ of $G$. Write $\Gamma = \text{Gal}(K_\infty/K)$ for the Galois group of the cyclotomic $\mathbb{Z}_p$ extension. We can and will identify $\Gamma$ with the Galois group of the cyclotomic $\mathbb{Z}_p$-extension of $\mathbb{Q}$. Let $\gamma$ denote a fixed topological generator of $\Gamma$. Put $\mathcal{O} = \mathbb{Z}_p[\chi]$ and $\Lambda = \mathcal{O}[[\Gamma]] = \mathcal{O}[[T]]$. For a finite order character $\rho : \Gamma \to \mathbb{C}^\times$, we define $\zeta \in \mu_{p^\infty}$ by $\zeta = \rho(\gamma)$.

---

[1] In Mazur's original situation, the proof is concluded (page 142-143) by an application of the Atkin-Lehner involution $w$, which interchanges the two components. This is not applicable here, as $g$ is an oldform, and hence will not in general be an eigenform for $w_M$.

Now let $f = \sum a_n q^n$ denote a weight-two cuspform for the group $\Gamma_1(N)$. We assume that $f$ is a simultaneous eigenform for all the Hecke operators of level $N$, and that $a_1 = 1$. We do *not* assume that $f$ is a newform. With these notations, the $\chi$-twisted $p$-adic L-function of $f$ is defined to be a power series

$$\mathscr{L}(f, \chi, T) \in \Lambda \otimes \mathbb{Q}_p$$

satisfying the following interpolation property for every nontrivial character $\rho$ of $\Gamma$:

$$\mathscr{L}(f, \chi, \zeta - 1) = \tau(\chi^{-1}\rho^{-1}) \cdot \alpha_p(f)^{-m} \cdot \frac{L(f, \chi\rho, 1)}{(-2\pi i)\Omega_f^\alpha}. \tag{24}$$

Here $\tau(\chi^{-1}\rho^{-1})$ denotes the usual Gauss sum attached to $\chi^{-1}\rho^{-1}$, and $p^m$ is the conductor of $\rho$. The quantity $\alpha_p(f)$ is the eigenvalue of $U_p$ on the $p$-stabilized newform associated to $f$. The period $\Omega_f^\alpha$ is a nonzero complex number, which we regard as fixed. The sign $\alpha$ is determined by $\pm 1 = \chi(-1)$. This interpolation property characterizes $\mathscr{L}(f, \chi, T)$, by the Weierstrass preparation theorem. If $\Sigma_0$ is any finite set of primes with $p \notin \Sigma_0$, then we define a nonprimitive L-function $\mathscr{L}^{\Sigma_0}(f, \chi, T)$, characterized by the interpolation formula

$$\mathscr{L}^{\Sigma_0}(f, \chi, \zeta - 1) = \tau(\chi^{-1}\rho^{-1}) \cdot \alpha_p(f)^{-m} \cdot \frac{L_{\Sigma_0}(f, \chi\rho, 1)}{(-2\pi i)\Omega_f^\alpha}, \tag{25}$$

where the $L_{\Sigma_0}(f, \chi\rho, s)$ denotes the complex L-function of $f$, stripped of the Euler factors at primes contained in $\Sigma_0$. The $p$-adic L-function of an elliptic curve is defined in a similar manner. Namely, one takes the $p$-adic L-function of the corresponding modular form $f$, and specifies the period by replacing $(-2\pi i)\Omega_f^\alpha$ with the Néron period $\Omega_E^\alpha$. Thus the L-functions of isogenous elliptic curves will differ by a constant.

**Proposition (3.7)** *Assume either that $E$ is optimal, or that $E[p]$ is irreducible. Assume that the character $\chi$ is unramified at all primes dividing $N$, and that $\chi$ is tamely ramified at $p$. Then the L-function $\mathscr{L}(E/\mathbb{Q}, \chi, T)$ is integral, i.e., we have $\mathscr{L}(E/\mathbb{Q}, \chi, T) \in \Lambda$.*

*Proof.* Let $f$ denote the newform associated to $E$. We see from (3.1) and (3.5) above that the numbers $\Omega_E^\alpha$ are such that the cocycle $\delta^\alpha$ defined by $\Omega_E^\alpha \delta^\alpha = \omega_f^\alpha$ lies in the integral cohomology group $H^1(\Gamma, \mathbb{Z}_p)$ (see (19); we are *not* assuming that $\alpha$ is admissible). In view of the hypotheses on $\chi$, one can check that, if $\zeta \neq 1$, then the values $\mathscr{L}(E, \chi, \zeta - 1)$ are given by the cap product of the integral class $\delta^\alpha$ against an appropriate *integral* element in the homology group $H_1(X_1; \mathbb{Z}[\chi, \zeta])$. This is well-known if the level $N$ is prime to $p$ (see [Ste82], Remark 1.6.2). If $p$ divides $N$, then $N$ must be divisible by precisely the first power of $p$, as $E$ has ordinary reduction. In this case a very similar argument to that of Stevens already cited proves the required integrality, when $\zeta \neq 1$. We omit the details. It follows that $\mathscr{L}(E, \chi, \zeta - 1)$ is $p$-integral for all $\zeta \neq 1, \zeta \in \mu_{p^\infty}$. The required result now follows from the Weierstrass preparation theorem. ∎

It is widely believed that the $p$-adic L-function of an elliptic curve is *always* represented by an integral power series. The defect in the foregoing proposition is that it gives no information about non-optimal curves $E$ admitting rational $p$-isogenies. However, we can remedy this defect if $E$ satisfies the hypotheses of theorem (1.3). Somewhat more generally, we have

**Corollary (3.8)** *Assume that $E$ admits a cyclic $p$-isogeny with kernel $\Phi$, such that $\Phi$ is either unramified at $p$ and odd, or ramified at $p$ and even, as a Galois module. Then the $\chi$-twisted $p$-adic L-function of $E$ is represented by an integral power series, for any even character $\chi$. If $E^{opt}$ is the optimal curve in the isogeny class of $E$, then the period $\Omega_E^+$ coincides with $\Omega_{E^{opt}}^+$, up to $p$-adic unit.*

*Proof.* It suffices to verify the final assertion. But this follows by exactly the same argument as was used in the conclusion of the proof of Lemma (3.6). Namely, one shows that the periods in question coincide with that of the minimal curve in the isogeny class. ∎

**Remark (3.9)** The foregoing proposition and corollary may be reformulated as follows. We want to show that the $\chi$-twisted L-function of $E$ is always represented by an integral power series. What we have shown is that this is in fact the case, *except* possibly for non-optimal curves $E$, and characters $\chi$ such that the sign determined by $\chi(-1)$ is not admissible for $E$. It is easy to see that, even in this latter case, the $p$-adic L-function will be integral if $E$ occurs as a subvariety of the Jacobian of some modular curve of level $N$. This includes for example the strong Weil curve, which in general will not be optimal. A complete resolution of the integrality question would follow from a conjecture of Stevens, [Ste89]. Namely, it would suffice to know that, if $E$ is any elliptic curve over $\mathbb{Q}$, then there exists an étale isogeny $E_{\mathbb{Z}_p}^{\text{opt}} \to E_{\mathbb{Z}_p}$, for the optimal curve $E^{\text{opt}}$ isogenous to $E$.

Our next result is the analytic ingredient in the proof of (1.4).

**Theorem (3.10)** *Let $E_1$ and $E_2$ be elliptic curves of level $N_1$ and $N_2$ respectively, such that Galois modules $E_1[p]$ and $E_2[p]$ are irreducible and isomorphic. Assume that the $E_i$ have good ordinary or multiplicative reduction at $p$. Let $\Sigma_0$ be the set of primes $q \neq p$ such that $q|N_1N_2$. Then there exists an element $u \in \mathcal{O}^\times$ such that we have the congruence $\mathscr{L}^{\Sigma_0}(E_1/\mathbb{Q}, \chi, T) \equiv u \cdot \mathscr{L}^{\Sigma_0}(E_2/\mathbb{Q}, \chi, T) \pmod{\pi\Lambda}$, for every character $\chi$, where $\pi$ is a uniformizing element in $\mathcal{O} = \mathbb{Z}_p[\chi]$.*

*Proof.* Let $M = \text{l.c.m}(N_1, N_2)$. Let $f_1 = \sum a_n q^n$ and $f_2 = \sum b_n q^n$ denote the modular forms associated to $E_1$ and $E_2$ respectively. Then the hypothesis that $E_1[p] \cong E_2[p]$ implies that we have the congruence $a_n \equiv b_n \pmod{p}$, whenever $(n, M) = 1$. This is clear if $(n, Mp) = 1$, and if $n = p$, then it follows from the results of [Wil88], as the curves $E_i$ are assumed to be ordinary. Now let $g_1 = \sum a'_n q^n$ and $g_2 = \sum b'_n q^n$ denote the eigenforms obtained from $f_1$ and $f_2$ respectively by dropping all primes in $\Sigma$. Then we have $a'_n \equiv b'_n \pmod{p}$, for *all* integers $n$. Theorem 1.10 of [Vat97] now yields a congruence as in the theorem, but the periods appearing will be the canonical periods attached to $g_1$ and $g_2$. The result follows from Lemma (3.6). ∎

Now we want to prove a similar theorem relating the $p$-adic L-function for an elliptic curve which admits a $\mathbb{Q}$-isogeny of degree $p$ to the $p$-adic L-function of a certain Eisenstein series, or equivalently, to the product of certain Kubota-Leopoldt $p$-adic L-functions. More precisely, we want to prove the congruence (12) stated in the introduction. Recall that $C = \mu_{p^\infty} \otimes \psi^{-1}$ and $D = (\mathbb{Q}_p/\mathbb{Z}_p) \otimes \psi$, where $\psi$ is an odd character with values in $\mathbb{Z}_p^\times$. (For our application, we assume $\psi$ is unramified at $p$.) Just as in theorem (3.10), we will prove a more general result by allowing a twist by a Dirichlet character $\chi$. We will assume that $\chi$ is even. The $p$-adic L-function $\mathscr{L}(C, \chi, T) \in \Lambda$ is characterized by the interpolation property

$$\mathscr{L}(C, \chi, \zeta - 1) = L(C, \chi\rho, 1) = L(\chi\psi^{-1}\rho, 0) \qquad (26)$$

for every nontrivial character $\rho$ of $\Gamma = \text{Gal}(\mathbb{Q}_\infty/\mathbb{Q})$. As before, $\zeta = \rho(\gamma)$ where $\gamma$ is a fixed topological generator of $\Gamma$. $\mathscr{L}(C, \chi, T)$ is related to the Kubota-Leopoldt $p$-adic L-function $L_p(\chi\omega\psi^{-1}, s)$ by

$$L_p(\chi\omega\psi^{-1}, s) = \mathscr{L}(C, \chi, \kappa(\gamma)^{-s} - 1)$$

for all $s \in \mathbb{Z}_p$. (Recall that $\kappa(\gamma) \in 1 + p\mathbb{Z}_p$ gives the action of $\gamma$ on $\mu_{p^\infty}$ when we identify $\Gamma$ with $\text{Gal}(\mathbb{Q}(\mu_{p^\infty})/\mathbb{Q}(\mu_p))$. Also, note that $\omega\psi^{-1} = \varphi$.) The Ferrero-Washington theorem asserts that $\mathscr{L}(C, \chi, T) \notin p\Lambda$ and the Mazur-Wiles theorem implies that the $\lambda$-invariant of $\mathscr{L}(C, \chi, T)$ is equal to $\text{corank}_\mathcal{O}(S_{C\otimes\chi}(\mathbb{Q}_\infty))$, which we denoted by $\lambda_{\chi\omega\psi^{-1}}$ in section 1. The notation $C \otimes \chi$ refers to the $\mathcal{O}$-module $C \otimes_{\mathbb{Z}_p} \mathcal{O}(\chi)$, where $G_\mathbb{Q}$ acts on the second factor by $\chi$. To obtain the nonprimitive $p$-adic L-function $\mathscr{L}^{\Sigma_0}(C,\chi,T)$, one multiplies $\mathscr{L}(C,\chi,T)$ by the $l$-th Euler factors $1 - \chi\psi^{-1}(l)(1+T)^{f_l}$ for each $l \in \Sigma_0$. Here $f_l \in \mathbb{Z}_p$ is determined by $\gamma_l = \gamma^{f_l}$, where $\gamma_l$ is the Frobenius automorphism for $l$ in $\Gamma$. The value at $T = \zeta - 1$ is the $l$-th Euler factor $1 - \chi\psi^{-1}\rho(l)$ in $L(\chi\psi^{-1}\rho, s)$ at $s = 0$. The $\lambda$-invariant of $\mathscr{L}^{\Sigma_0}(C,\chi,T)$ is $\lambda_{\chi\omega\psi^{-1},\Sigma_0} = \lambda_{\chi\varphi,\Sigma_0}$.

The $p$-adic L-function $\mathscr{L}(D,\chi,T) \in \Lambda$ is characterized by the interpolation property

$$\begin{aligned}\mathscr{L}(D,\chi,\zeta-1) &= \tau(\chi^{-1}\psi^{-1}\rho^{-1})L(D,\chi\rho,1)/2\pi i \\ &= \tau(\chi^{-1}\psi^{-1}\rho^{-1})L(\chi\psi\rho,1)/2\pi i \\ &= \tfrac{1}{2}L(\chi^{-1}\psi^{-1}\rho^{-1},0)\end{aligned} \qquad (27)$$

for every nonzero character $\rho$ of $\Gamma$, where $p^m$ is the conductor of $\rho$. The last equality follows from the functional equation. $\mathscr{L}(D,\chi,T)$ is related to the Kubota-Leopoldt $p$-adic L-function $L_p(\omega\chi^{-1}\psi^{-1},s)$ by

$$L_p(\omega\chi^{-1}\psi^{-1},s) = \tfrac{1}{2}\mathscr{L}(D,\chi,\kappa(\gamma)^s - 1)$$

for all $s \in \mathbb{Z}_p$. The $\mu$-invariant of $\mathscr{L}(D,\chi,T)$ is again zero and its $\lambda$-invariant is $\lambda_{\omega\chi^{-1}\psi^{-1}} = \lambda_{\chi\psi}$, which is equal to $\text{corank}_\mathcal{O}(S_{D\otimes\chi}(\mathbb{Q}_\infty))$. To obtain $\mathscr{L}^{\Sigma_0}(D,\chi,T)$, one multiplies by the Euler factors $1 - \chi\psi(l)l^{-1}(1+T)^{f_l}$ for all $l \in \Sigma_0$. The value at $T = \zeta - 1$ is the $l$-th Euler factor $1 - \chi\psi\rho(l)l^{-1}$ in $L(\chi\psi\rho,s)$ at $s = 1$. The $\lambda$-invariant of $\mathscr{L}^{\Sigma_0}(D,\chi,T)$ is $\lambda_{\chi\psi,\Sigma_0}$.

Suppose that $\varphi$ and $\psi$ are as before, namely the $\mathbb{Z}_p^\times$-valued characters corresponding to the composition factors $\Phi$ and $\Psi$ in the $G_\mathbb{Q}$-module $E[p]$, assuming that $E$ admits a $\mathbb{Q}$-isogeny of degree $p$. We assume that $\psi$ is odd and unramified at $p$, or equivalently that $\varphi$ is even and ramified at $p$. In this case, the admissible sign is plus, and the canonical period is the real period of $E$ (up to multiplication by a $p$-adic unit). Let $N$ denote the level of $E$ and let $\Sigma_0$ denote any set of primes containing all primes $l \neq p$ dividing $N$, but not including $p$. Define $G = \Sigma b_n q^n$ to be the weight-two Eisenstein series determined by $\Sigma b_n n^{-s} = L^{\Sigma_0}(\psi,s)L^{\Sigma_0}(\psi^{-1},s-1)$. Here the superscript $\Sigma_0$ indicates that we consider the non-primitive L-functions obtained by omitting the Euler factors for all $l \in \Sigma_0$.

For each even Dirichlet character $\psi$, we let $\mathscr{L}(G,\chi,T)$ denote the $p$-adic L-function (associated to $G$ and $\chi$) characterized by the interpolation property

$$\begin{aligned}\mathscr{L}(G,\chi,\zeta-1) &= \tau(\chi^{-1}\psi^{-1}\rho^{-1})L(G,\chi\rho,1)/2\pi i \\ &= (L^{\Sigma_0}(\chi\psi^{-1}\rho,0))(\tau(\chi^{-1}\psi^{-1}\rho^{-1})L^{\Sigma_0}(\chi\psi\rho,1)/2\pi i)\end{aligned}$$

for all nontrivial characters $\rho$ of $\Gamma$, where as before $\rho(\gamma) = \zeta$. Then we clearly have

$$\mathscr{L}(G,\chi,T) = \mathscr{L}^{\Sigma_0}(C,\chi,T)\mathscr{L}^{\Sigma_0}(D,\chi,T). \qquad (28)$$

The $\mu$-invariant of $\mathscr{L}(G,\chi,T)$ is zero because that is true for each factor in (28). Also, the $\lambda$-invariant of $\mathscr{L}(G,\chi,T)$ is equal to $\lambda_{\chi\varphi,\Sigma_0} + \lambda_{\chi\psi,\Sigma_0}$. The two terms are the $\mathcal{O}$-coranks of $S_{C\otimes\chi}^{\Sigma_0}(\mathbb{Q}_\infty)$ and $S_{D\otimes\chi}^{\Sigma_0}(\mathbb{Q}_\infty)$, respectively. Theorem (1.3) is a consequence of the congruence in the following theorem. We just take $\chi$ to be the trivial character. We then obtain that $\lambda_{E,\Sigma_0}^{\text{anal}} = \lambda_{\varphi,\Sigma_0} + \lambda_{\psi,\Sigma_0}$ which is in turn equal to $\lambda_{E,\Sigma_0}^{\text{alg}}$ by (16). Thus $\lambda_{E,\Sigma_0}^{\text{alg}} = \lambda_{E,\Sigma_0}^{\text{anal}}$ which by theorem (1.5) implies that $\lambda_E^{\text{alg}} = \lambda_E^{\text{anal}}$. The vanishing of $\mu_E^{\text{anal}}$ also follows from the congruence. (The vanishing of $\mu_E^{\text{alg}}$ was proved in section 2 under the hypotheses of theorem (1.3).)

**Theorem (3.11)** *Let $\chi$ be any even character. Then we have congruence*

$$\mathscr{L}^{\Sigma_0}(E/\mathbb{Q},\chi,T) \equiv u\mathscr{L}(G,\chi,T) \pmod{\pi\Lambda},$$

*where $u$ is a unit in $\mathcal{O}$.*

*Proof.* Let $f$ denote the cuspform associated to $E$, and let $g$ denote the form obtained by removing all Euler factors at primes $q \in \Sigma$. Then both $f$ and $G$ are simultaneous eigenforms at some common level $M$. Furthermore, if $g = \sum a_n q^n$, then we we have the congruence $a_n \equiv b_n \pmod{p}$ for every integer $n$. We contend now that the constant term of $G$ at every cusp of $X_1(M)$ is divisible by $p$. Observe first that the constant term $b_0$ of $G$ at infinity vanishes; this follows from the fact that $L(s,G)$ is holomorphic ($\chi \neq 1$), and the well-known characterization of this constant term as the residue of the L-function at $s = 1$. Thus the modular form $g - G = \sum_{n\geq 1} c_n q^n$ is such that $p | c_n$ for all $n$. If $(p,M) = 1$, then our contention follows immediately from the $q$-expansion principle. If $(p,M) \neq 1$, then $M$ must be divisible by precisely the first power of $p$. In this case one can argue as follows. The $q$-expansion principle ensures only that $g - G$ vanishes on the component containing infinity of $X_1(N)_k$, where $k$ is the residue field of $\mathbb{Z}_p[\zeta_N]$, and $\zeta_N = e^{2\pi i/N}$ is a primitive $N$-th root of unity. [2] This component contains the images of the so-called "infinity cusps," which are those represented by rational numbers of the form $a/p^r$, with $r \geq 1$ and $a$ prime to $p$. It follows that the Eisenstein series $G$ has the property that its constant term vanishes modulo $p$ at the infinity cusps. The assertion about the constant terms now follows from Hida's determination of the ordinary Eisenstein series ([Hid85], Thm. 5.8), since $g$ is obviously ordinary. We have therefore checked all the hypotheses in Theorem (2.10) of [Vat97]. Applying that theorem gives the congruence

$$\mathscr{L}(g,\chi,\zeta - 1) \equiv \mathscr{L}(G,\chi,\zeta - 1) \pmod{\pi},$$

for every $\zeta \neq 1$. Our theorem follows as before. ∎

---

[2] Here we are using the conventions of [MW84], Ch. 2; note that the curve denoted there by $X_1(N)_\mathbb{Q}$ is not the same as the one considered in 3. However, the two become isomorphic over $\mathbb{Q}(\zeta_N)$.

## Modular Forms with non-rational coefficients

Finally, we want to briefly discuss the situation for $p$-adic L-functions associated to modular forms whose Fourier coefficients are not necessarily rational numbers. Much of the theory extends to this context, but there are some important differences, which we will now describe. Thus let $f = \sum a_n q^n$ denote a weight 2 eigenform on $\Gamma_1(N)$, with coefficients in the $p$-adic integer ring $\mathcal{O}$. We write $\rho_f$ for the usual Galois representation $G_\mathbb{Q} \to GL_2(\mathcal{O})$ attached to $f$, so that $\operatorname{Tr}(\operatorname{Frob}(q)) = a_q$, for all primes $q$ not dividing $Np$. We assume as usual that $N$ is divisible by at most the first power of $p$, and that $a_p$ is a unit in $\mathcal{O}$. We have seen in the preceding sections how to define a Selmer group $\operatorname{Sel}_f(\mathbb{Q}_\infty)$ for $\rho_f$, and how to define a $p$-adic L-function $\mathscr{L}(f, T)$. Then the main conjecture states that $\operatorname{Sel}_f(\mathbb{Q}_\infty)$ is $\Lambda$-cotorsion, and has characteristic ideal generated by the power series $\mathscr{L}(f, \mathbf{T})$. We can define Iwasawa invariants $\mu_f^{\text{alg}}, \mu_f^{\text{anal}}, \lambda_f^{\text{alg}}$, and $\lambda_f^{\text{anal}}$ as before, just as in the case of elliptic curves. Our task is, once again, to prove that $\mu_f^{\text{anal}} = \mu_f^{\text{alg}} = 0$, and that $\lambda_f^{\text{anal}} = \lambda_f^{\text{alg}}$. We emphasize here that the definition of the $p$-adic L-function involves the choice of a complex period for $f$. From the viewpoint of modular forms, the natural choice is the canonical period of [Vat97]. However, the form $f$ is associated to an abelian variety $A$, and one can also define a period in terms of the geometry of $A$ (we will describe this below) although the definition is not so canonical as in the case of elliptic curves. These definitions can be shown to be equivalent in most cases; presumably the periods are equal in general, but we cannot prove this at present.

Consider first the analogue of Theorem (1.3). Let $\mathbf{T}$ denote the Hecke ring for $\Gamma_1(N)$, with coefficients in $\mathcal{O}$, and let $\mathfrak{m}$ denote the maximal ideal determined by $f$, and let $\mathbb{F}$ denote the residue field of the completion $\mathbf{T}_\mathfrak{m}$. There is a semisimple representation $\rho_m : G_\mathbb{Q} \to GL_2(\mathbb{F})$ satisfying $\operatorname{Tr}(\operatorname{Frob}(q)) = T_q$, for any $q$ with $(Np, q) = 1$. The analogue to the hypothesis of Theorem (1.3) is that $\rho_\mathfrak{m}$ be *reducible* in the sense that there exist characters $\varphi, \psi : G_\mathbb{Q} \to \mathbb{F}^\times$, such that

$$\rho_\mathfrak{m} = \varphi \oplus \psi.$$

We would like to define canonical periods for $f$, but unfortunately, this is not always possible. We will therefore make the following assumption: *the characters $\varphi$ and $\psi$ are distinct when restricted to the decomposition group $D_p$*. Such a representation is said to be $p$-distinguished (see Theorem 2.1 in [Wil95], Chapter 2; also [Vat97], Theorem 2.7.) This condition will always be satisfied if, for instance, the level $N$ is prime to $p$. With the assumption that $\rho_\mathfrak{m}$ is $p$-distinguished, we single out the character $\psi$ by requiring that it be unramified, and that it satisfy $\psi(\operatorname{Frob}(p)) = a_p$ in $\mathbb{F}^\times$. Let $\alpha = \pm$ be the choice of sign determined by $-\psi(-1) = \pm 1$. Since $\rho_\mathfrak{m}$ is assumed to be $p$-distinguished, Theorem 2.7 of [Vat97] implies that the canonical period $\Omega_f^\alpha$ exists.

Let $\tilde{\varphi}$ and $\tilde{\psi}$ denote the Teichmuller lifts of $\varphi$ and $\psi$ to $\mathcal{O}^\times$. There exists an Eisenstein series $G = \sum b_n q^n$ of level $N$ such that $b_n = \tilde{\psi}(q) + q\omega^{-1}(q)\tilde{\varphi}(q)$, for all primes $q$ with $(q, Np) = 1$, and $a_n \equiv b_n \pmod{\pi}$ for all $n > 0$. Furthermore, we may assume that the constant term of $G$ at every cusp is divisible by the uniformizer $\pi$ of $\mathcal{O}$. (This is a standard argument; see the proof of Theorem 3.3 in [Vat97], or the proof of Theorem (3.11) above.) Let $\chi$ be any Dirichlet character with conductor prime to $N$, and with $\chi(-1) = -\psi(-1)$. We may choose $\mathcal{O}$ large enough to contain the values of $\chi$. Let $\Sigma_0$ denote any finite set of primes, with $p \notin \Sigma_0$, and containing all other primes $q \neq p$ dividing $N$. As in the proof of Theorem (3.11), we obtain the following result:

**Theorem (3.12)** *Assume that $\rho_\mathfrak{m}$ is $p$-distinguished. Then there exists an invertible power series $U_\chi(T)$ such that the following congruence holds:*

$$\mathscr{L}^{\Sigma_0}(f, \chi, T) \equiv U_\chi(T) \cdot \mathscr{L}^{\Sigma_0}(\chi\varphi\omega^{-1}, T) \cdot \mathscr{L}^{\Sigma_0}(\chi^{-1}\psi^{-1}, (1+T)^{-1} - 1) \pmod{\pi\Lambda}.$$

*Thus, the invariant $\mu_f^{anal}$ is trivial. The main conjecture is true for $f$.*

Here the period appearing in $\mathscr{L}^{\Sigma_0}(f, \chi, T)$ is the canonical period $\Omega_g^\alpha$ attached to the eigenform $g$ obtained from $f$ by dropping all Euler factors in $\Sigma_0$. It is possible to give a geometric interpretation of this quantity as a period on a suitable abelian variety. To this, one can proceed as follows. Without loss of generality, we may assume that $f$ is a newform of level $N$. Let $A = A_f$ denote the abelian variety quotient of $J = J_1(N)$ constructed by Shimura. We may assume that $A$ is optimal in the sense that the map $J \to A$ has connected kernel. Let $R \subset \mathbb{C}$ denote the $\mathbb{Z}$-algebra generated by the Fourier coefficients of $f$. Then $R$ is an order is a number field, and there is an embedding of $R$ into the endomorphism ring of $A$. Working locally at $p$, it can be shown that the cotangent space $\text{Cot}(J) = H^0(X_1, \Omega^1)$ is locally free over $R_p = R \otimes \mathbb{Z}_p$. Furthermore, each of the spaces $H_1(J, \mathbb{Z}_p)^\pm$ is also free over $R_p$. Let $\omega$ denote a generator for $\text{Cot}(J) \otimes \mathbb{Z}_p$ as an $R_p$ module, and let $c^\pm$ denote a generator for $H_1(J, \mathbb{Z}_p)^\pm$. Then we can define geometric periods $\Omega_A^\pm$ by $\Omega_A^\pm = \int_{c^\pm} \omega$. If the representation $\rho_\mathfrak{m}$ is $p$-distinguished and if the abelian variety has good or semistable reduction at $p$, it can be shown that the periods $\Omega_A^\alpha = \Omega_f^\alpha$, for an admissible sign $\alpha$. We note however that the abelian varieties that arise for forms on $J_1(N)$ need not have good or multiplicative reduction when $p|N$, even if $\rho_\mathfrak{m}$ is $p$-distinguished.

There is also an analogue of Theorem (1.4) for forms with non-rational Fourier coefficients. To formulate this, consider a pair of newforms $f = \sum a_n q^n$, $g = \sum b_n q^n$ of level $N$ and $M$ respectively. We assume that the ring $\mathcal{O}$ is sufficiently large to contain the Fourier coefficients of $f$ and $g$, and that both $a_p$ and $b_p$ are $p$-adic units. Suppose there exists a finite set of primes $\Sigma_0$, containing all primes $q \neq p$ dividing $NM$, but not containing $p$, such that $a_n \equiv b_n$ (mod $\pi$) for all $n$ indivisible by the primes in $\Sigma_0$. There exists a unique representation $\overline{\rho} : G_{\mathbb{Q}} :\to GL_2(\mathcal{O}/\pi\mathcal{O})$ satisfying $\text{Tr}\,(\rho(\text{Frob}(q)) = \overline{a}_q = \overline{b}_q$ in the residue field $\mathcal{O}/\pi\mathcal{O}$, for all primes $q \notin \Sigma_0$. Then $\overline{\rho}$ is ordinary by a theorem of Wiles (see [Wil88], Theorem 2.1.4). Let $A$ and $B$ denote the optimal abelian variety quotients of $J_1(N)$ and $J_1(M)$ associated by Shimura to $f$ and $g$ respectively. By proceeding as above, we may define geometric periods $\Omega_A^{\pm}$ and $\Omega_B^{\pm}$ for $A$ and $B$.

Let $\Sigma_0$ be a finite set of primes as before, containing all primes $q \neq p$ dividing $NM$, but not containing $p$. The results of section 2 in [Vat97] now imply the following result.

**Theorem (3.13)** *Let the hypotheses be as above. Assume in addition that the representation $\overline{\rho}$ is irreducible and p-distinguished. Then we have the congruence*

$$\mathscr{L}^{\Sigma_0}(f, \chi, T) \equiv \mathscr{L}^{\Sigma_0}(g, \chi, T) \pmod{\pi}.$$

*If the abelian variety $A$ has good or semistable reduction at p, then the geometric periods of $A$ coincide with the canonical periods of $f$ up-to p-adic unit. A similar statement holds for $B$ and $g$.*

Combining this with the the Selmer group calculations of Section 2 and the arguments outlined in the introduction, we obtain the following result.

**Corollary (3.14)** *Let the notation be as above. Assume that the representation $\overline{\rho}$ is p-distinguished and irreducible. If the equalities $\mu_f^{alg} = \mu_f^{anal} = 0$ and $\lambda_f^{alg} = \lambda_f^{anal}$ hold, then we have the further equalities $\mu_g^{alg} = \mu_g^{anal} = 0$ and $\lambda_g^{alg} = \lambda_g^{anal}$. The main conjecture holds for $g$.*